\documentclass[12pt]{article}
\usepackage{amssymb,graphicx,color}

\def\d{{\rm d}}
\def\mi{{\rm i}}
\def\eps{\varepsilon}
\def\g{\gamma}
\def\G{\Gamma}
\def\l{\lambda}
\def\L{\Lambda}
\def\z{\zeta}
\def\Re{\mathop{\rm Re\,}\nolimits}
\def\Im{\mathop{\rm Im\,}\nolimits}
\def\e{\mathop{\rm e}\nolimits}
\def\Res{\mathop{\rm Res}\nolimits}
\def\hf{{\textstyle{1 \over 2}}}
\def\qt{{\textstyle{1 \over 4}}}
\def\tq{{\textstyle{3 \over 4}}}
\def\of{{\textstyle{1 \over 5}}}
\def\tf{{\textstyle{2 \over 5}}}

\def\defi{\stackrel{\rm def}{=}}
\def\si{\!\!\! &}
\def\se{& \!\!\!}
\def\ni{\noindent}

\newcommand{\beq}{\begin{equation}}
\newcommand{\eeq}{\end{equation}}
\newcommand{\bea}{\begin{eqnarray}}
\newcommand{\eea}{\end{eqnarray}}
\newcommand{\ttfr}[2]{\frac{\mbox{\small $#1$}}{\mbox{\small $#2 \strut$}}}

\title{Discretized Keiper/Li approach\\
to the Riemann Hypothesis}

\author{{\bf Andr\'e Voros}\\
Institut de physique th\'eorique, Universit\'e Paris-Saclay, CEA, CNRS\\
F-91191 Gif-sur-Yvette Cedex (France)\\
E-mail: \tt andre.voros@ipht.fr}

\begin{document}

\maketitle

\begin{abstract}
The Keiper--Li sequence $\{ \lambda _n \}$ is most sensitive to the Riemann Hypothesis
asymptotically ($n \to \infty$), but highly elusive both analytically and numerically.
We deform it to fully explicit sequences, simpler to analyze and to compute (up to $n=5 \cdot 10^5$ by G. Misguich).
This also works on the Davenport--Heilbronn counterexamples, 
thus we can demonstrate explicit tests that selectively react to zeros \emph{off} the critical line. 
\end{abstract}

The present text develops our computations announced from 2015 \cite{VU} 
on a novel type of \emph{explicit} sequential criteria for the \emph{Riemann Hypothesis}.

The Riemann zeta function $\z(x) :$ e.g., \cite{Ti}\cite[chaps. 8, 15]{Da},
has countably many zeros $\rho $ in the upper half-plane $\{ \Im \rho >0 \}$,
which accumulate to $+ \mi \infty$ (and their complex-conjugates, which we will not count).
The Riemann Hypothesis (RH) \cite{Ri} states that \emph{all those zeros} lie on the \emph{critical line} 
${ \{ \Re \rho = \hf \} }$, and this remains a major open conjecture. 
 
Within \emph{bounded heights} $T= \Im \rho $ on the other hand, 
direct tests have always validated $\Re \rho \equiv \hf $, 
now up to the $10^{13}$-th zero \cite{G}. If $T_0$ denotes the highest $T$ up to
which RH is verified (at a given time), that roughly sets
\beq
\label{TV}
T_0 \approx 2.4 \cdot 10^{12} \quad \mbox{currently (as of 2004)} .
\eeq

We will examine how novel sequential criteria might probe RH \emph{further}, 
i.e., sense any zeros $\rho ' = \hf + t + \mi T$ with $t>0$ (violating RH) and $T>T_0$.

First (\S~1) we review the usual \emph{Keiper--Li} constants $\{ \l _n \}_{n \ge 1}$. 
Their (real-valued) sequence is RH-sensitive and could detect such $\rho '$
for ${n \gtrsim T^2/t}$ (currently: above $\approx 10^{25}$). 
But the $\l _n$ involve derivatives of $\log \z $ to any orders, which are most elusive 
both analytically and numerically (in which case the difficulty steeply rises with $n :$
only recently was $n=10^5$ reached). 

We propose instead (\S~2) to modify the Keiper--Li definition, by discretizing the above derivatives 
to \emph{finite differences} while retaining RH-sensitivity.
We thus obtain a variant sequence $\{ \L _n \}$ in the \emph{elementary explicit form} (\ref{EKL}), 
also more expedient to encode numerically (G. Misguich reached $n=5 \cdot 10^5$).

In \S~3 we confirm analytically that the sequence $\{ \L _n \}$ is RH-sensitive, 
in its asymptotic tail like $\{ \l _n \}$, just a bit less simply; 
we also extend the framework to certain Dirichlet L-series.

In \S~4 we initiate quantitative and numerical studies.
We find two main thresholds for the visibility of the above $\rho '$ through $\{ \L _n \} :$ 
for $n \gtrsim T^{1+2/t}$ (currently: above $\approx 10^{60}$) detection is \emph{certain};
inversely for $n \lesssim \hf T \e^{1/t}$ (currently $\approx 10^{13}$ or higher) 
detection is outright \emph{barred} by an uncertainty principle. 
Somewhere in-between, achievable sensitivity will depend on how deeply
$\{ \L _n \}$ can be described and processed in the future.
In any case and like $\{ \l _n \}$, $\{ \L _n \}$ probes RH
in a \emph{complementary} manner to the classic direct methods (Riemann--Siegel based),
currently more efficient but strictly local in~$T$.
Thanks to the Davenport--Heilbronn counterexamples, 
we can concretely demonstrate how a (generalized) $\L _n $-sequence signals $\Re \rho ' \ne \hf $.

We finally analyze in detail the loss of precision in the numerical calculation of $\L _n$.
It grows linearly with $n$, being on average twice that for $\l _n $, 
but for $\L _n $ that is the \emph{sole} computational stumbling block, and a \emph{purely logistic} one 
(arbitrary-precision programs are fully ready, only computing power has to rise to the task).
Hence overall, $\{ \L _n \}$ is easier to calculate, has a lower uncertainty-principle barrier,
and its study only begins. This discretized version of the Keiper--Li sequence thus ought to be
more efficient in the long run than the original $\{ \l _n \}$ for testing RH further.

An Appendix describes a similar discretization for a variant (more symmetrical) Keiper--Li sequence.

\section{Main notations and background}
\label{MNB}

\ni $\bullet \ 2 \xi (\cdot ):$ a completed zeta function, obeying Riemann's Functional Equation,
\beq
\label{CZ}
2 \xi (x) \defi x(x-1) \pi^{-x/2} \, \G (x/2) \, \z (x) \equiv 2 \xi (1-x) :
\eeq
$2 \xi $ is better normalized for us than Riemann's $\xi $-function: $2 \xi (0) = 2 \xi (1) = 1$.

\ni The symmetry in (\ref{CZ}) makes us also denote $x = \hf +t +\mi T \ $ ($t,T$ real).
\smallskip

\ni $\bullet \ \{\rho \}:$ the set of zeros of $\xi $ 
(the \emph{nontrivial zeros} of $\z $ or \emph{Riemann zeros}, counted with multiplicities if any);
this set has the symmetry axes $\mathbb R$ and ${L \defi \{ \Re x = \hf \}}$ (the \emph{critical line}),
it lies in the \emph{critical strip} ${ \{ 0< \Re x <1 \} }$,
and its ordinate-counting function $N(T) \defi \# \{ \rho \mid 0 \le \Im \rho \le T \}$ 
has the \emph{Riemann--von Mangoldt} asymptotic form: 
\beq
\label{RvM}
N(T) = \overline N (T) + O(\log T) \quad \mbox{for } T \to +\infty , \quad
\mbox{with } \overline N (T) = \frac{T}{2\pi} \log \frac{T}{2\pi {\rm e}} \, ;
\eeq
infinite summations over zeros are meant to be ordered symmetrically, as
\beq
\label{SG}
\sum_\rho \defi \lim_{T \to \infty} \sum_{|\Im \rho| \le T} .
\eeq
\smallskip

We will dub $\rho '$ any zeros for which we assume $\Re \rho ' > \hf $ (violating RH).

\smallskip

\ni $\bullet \ \z (x,w) \defi \sum\limits_{k=0}^\infty (k+w)^{-x} :$ the Hurwitz zeta function.

\ni $\bullet \ \G (\cdot) :$ the Euler Gamma function, \quad $\g :$ Euler's constant ($\approx 0.5772156649$).

\ni $\bullet \ k!! :$ the double factorial, used here for \emph{odd} integers $k$ only, in which case
\bea
k!! \si\defi\se k(k-2) \cdots 1 \qquad \qquad \quad \mbox{for odd } k>0 , \nonumber \\
\si\defi\se 2^{(k+1)/2} \, \G (\hf k+1) / \sqrt \pi \quad \mbox{for odd } k\gtrless 0 
\quad (e.g.,\ (-1)!!=1, \ (-3)!!=-1) . \nonumber 
\eea
\ni $\bullet \ (j,k,\ell )! \defi (j+k+\ell )! \, /\, (j! \, k! \, \ell !)$ \ (3-index multinomial coefficients).

\ni $\bullet \ B_{2m} :$ Bernoulli numbers,  $B_{2m} (\cdot) :$ Bernoulli polynomials \ ($m=0,1,2,\ldots $).

\subsection{The Keiper and Li coefficients}

In 1992 Keiper \cite{K} considered a real sequence $\{ \l _n \}$ of generating function
\beq
\label{KEd}
\varphi (z) \defi \log 2 \xi (M(z)) \equiv \sum_{n=1}^\infty \l ^{\rm K}_n \, z^n ,
\qquad M(z) \defi \frac{1}{1-z}, 
\eeq
(we write $\l ^{\rm K}_n$ for \emph{Keiper's} $\l _n$), deduced that
\beq
\label{KEp}
\l ^{\rm K}_n \equiv n^{-1} \sum_\rho \bigl[ 1-(1-1/\rho)^n \bigr] ,
\eeq
and that RH $\Rightarrow \ \l ^{\rm K}_n >0 \ (\forall n)$, then wrote (without proof nor elaboration):
``In fact, if we assume the Riemann hypothesis, and further that the zeros are very evenly distributed, 
we can show that
\beq
\label{Kap}
\l _m \approx \frac{\log m}{2} -\frac{\log (2\pi ) +1-\g }{2} . \, {" \atop } \quad
\mbox{[for $\l ^{\rm K}_m \,$; \ $m \to \infty$ seems implied.}]
\eeq
In (\ref{KEd}), the conformal mapping $M : z \mapsto x$ acts to pull back each zero $\rho $ to
\beq
z_\rho \defi M^{-1} (\rho ) = 1-1/\rho ,
\eeq
and \emph{the critical line $L$ to the unit circle} $\{ |z|=1 \}$, ensuring that (Fig.~1):
\beq
\label{KLS}
\mbox{RH} \iff |z_\rho | \equiv 1 \ (\forall \rho ) \iff \varphi \mbox{ is \emph{regular in the whole disk} } \{ |z|<1 \}.
\eeq
\smallskip

In 1997 Li \cite{LI1} independently introduced another sequence $\{ \l _n \}$, through
\beq
\label{LId}
\l ^{\rm L}_n = \frac{1}{(n-1)!}\, \frac{\d ^n}{\d x^n} \bigl[ x^{n-1} \log 2 \xi (x) \bigr]_{x=1},
\quad n=1,2,\ldots ,
\eeq
(we write $\l ^{\rm L}_n$ for \emph{Li's} $\l _n$), 
and proved the \emph{sharp equivalence (Li's criterion)}: 
\beq
\label{LiC}
\mbox{RH} \quad \iff \quad \l ^{\rm L}_n \ge 0 \mbox{ for all }n.
\eeq
He also obtained $\l ^{\rm L}_n \equiv \sum\limits_\rho \bigl[ 1-(1-1/\rho)^n \bigr] $,
which in view of (\ref{KEp}) entails 
\beq
\l ^{\rm L}_n \equiv n \, \l ^{\rm K}_n \qquad \mbox{for all } n=1,2,\ldots ;
\eeq
our superscripts K vs L will disambiguate $\l _n$ whenever the factor $n$ matters.

\subsection{Probing RH through the Keiper--Li constants $\l _n$}
\label{PRH}

Li's criterion (\ref{LiC}) makes it clear that the Keiper--Li sequence is \emph{RH-sensitive}: 
how efficiently can it then serve to test RH?

Known results actually imply that, beyond the present frontier (\ref{TV}), 
the sequence $\{ \l _n \}$ may effectively probe RH 
only in its \emph{tail} $n \gg 1$ and via its \emph{asymptotic form} for $n \to \infty$.
\smallskip

In 2000 Oesterl\'e \cite[prop.~2]{O} proved (but left unpublished \cite[\S~2.3]{BPY}) that 
\beq
\label{OP}
\Re \rho =\hf \mbox{ for all zeros with } |\Im \rho | \le T_0 \ \ \Longrightarrow
\ \ \l _n \ge 0 \mbox{ for all } n \le T_0^{\, 2} ,
\eeq
and that under RH, \cite[\S~2]{O} 
\beq
\label{OE}
\l ^{\rm L}_n = n (\hf \log n + c) + o(n)_{n \to \infty} \quad 
\mbox{with } c = \hf (\g -\log 2\pi -1) ,
\eeq
which concurs with Keiper's formula (\ref{Kap}) but now assuming \emph{RH alone}.

In 2004 Ma\'slanka \cite{M1}\cite{M2} computed a few thousand $\l ^{\rm L}_n$-values numerically
and inferred asymptotic conjectures on them for the case RH true.

In 2004--2006, inspired by the latter (but unaware of \cite{O}), we used the saddle-point method 
to draw an \emph{asymptotic criterion} for RH \cite{V}: as $n \to \infty$,
\footnotetext[1]{Erratum: in \cite{V}\cite[chap.~11]{VB} 
we missed the overall ($-$) sign (with no effect on our conclusions), which we rectified in \cite{VK}.}
\bea
\label{ASF}
\bullet \ \mbox{RH false:} \si\se \l ^{\rm L}_n \sim - \!\sum_{\Re \rho ' > \frac{1}{2} } z_{\rho '}^{\, -n}
\pmod {o(r^{-n})\ \forall r<1} \\
\si\se \qquad \mbox{(exponentially growing oscillations \emph{with both signs});}\footnotemark \nonumber \\
\label{ASC}
\bullet \ \mbox{RH true: } \si\se \l ^{\rm L}_n \sim n (\hf \log n + c) \pmod {o(n)} \\
\si\se \qquad \mbox{(implying asymptotic \emph{positivity}).} \nonumber
\eea

Note: the remainder term $o(n)$ in (\ref{OE}) got improved to $O(\sqrt n \log n)$ 
by Lagarias \cite{Lg} (2007),
and to $n y_n$ with $\{y_n\} \in \ell^2$ by Arias de Reyna \cite{AR} (2011).
\smallskip

Then, either (\ref{OP}) or (\ref{ASF})--(\ref{ASC}) make searches through $\{ \l _n \}$  
for violations of RH increasingly harder as the floor height $T_0$ (currently (\ref{TV})) goes higher.

First, Li's sign test (\ref{LiC}) cannot, by (\ref{OP}), operate
below $n=T_0^{\, 2}$ to the least (${\approx 5 \cdot 10^{24}}$ today).
As for the asymptotic test (\ref{ASF})--(\ref{ASC}): if a zero $\rho ' = \hf +t+\mi T$ violates RH
(thus $|T |\ge T_0 \gg 1$), then its imprint $z_{\rho '}^{\, -n}$ in (\ref{ASF}) 
grows detectable against the background (\ref{ASC}) only for $n \gtrsim T^2 /|t|$ \cite[\S~11.3]{VB}.
Seen in the $z$-plane, that RH-violation is measured by $\delta |z| = |z_{\rho '}| -1$ (Fig.~1 right); 
now $\delta |z| / |z_{\rho '}| \approx -t/T^2$ as $|T| \gg 1$, 
therefore $n \gtrsim T^2 /|t|$ means $n \, |\delta \log |z_{\rho '}||\gtrsim 1$;
which is no less than the \emph{uncertainty principle},
as $(\mi \log z)$ and $n$ are Fourier-conjugate variables in (\ref{KEd})
(cf. also (\ref{LR}) below, where $\theta \equiv \mi ^{-1}\log z$;
simply, holomorphy as in $\varphi (z)$ implies $2\pi$-periodicity in this variable $\theta $,
and a power series in $z$ such as (\ref{KEd}) is but a Fourier series in~$\theta $). 
Now the uncertainty principle is universal, implying that
actually \emph{any} detection scheme of~$\rho '$ through $\l _n$ will require 
\beq
\label{Lup}
n \gtrsim T^2 /|t| .
\eeq
With $|t|<\hf $, the most favorable case ($t=\hf -0, \ T=T_0$) then necessitates
\beq
\label{LIt}
n \gtrsim 2T_0^{\, 2}, \quad \mbox{currently implying} \ \ n \gtrsim 10^{25} .
\eeq
It is then no wonder that published $\l _n$-plots (having $n \lesssim 7000$) 
solely reflect the RH-true pattern (\ref{ASC}) (already from $n \approx 30$) \cite{K}\cite{M1}.

So, whether one would take (\ref{LiC}) or (\ref{ASF})--(\ref{ASC}) to track violations of RH,
the Keiper--Li sequence $\{ \l _n \}$ only matters in its \emph{asymptotic tail} ${\{ n \gg 1 \}}$,
where the alternative (\ref{ASF})--(\ref{ASC}) rules (and enacts Li's sign property as well).

At the same time, the $\l _n$ are quite elusive analytically \cite{BL}\cite{C2}.
Numerically too (see Ma\'slanka \cite{M1}\cite{M2} and Coffey \cite{C1}),
their evaluation requires a recursive machinery, of intricacy blowing up with $n$;
\cite[fig.~6]{M2} moreover reports a loss of precision of $\approx 0.2$ decimal place per step~$n$ 
(when working \emph{ex nihilo} - i.e., using no Riemann zeros in input);
only $\l _n$-values up to $n \approx 4000$ were accessed that way
until recently, with $n=10^5$ attained by Johansson \cite[\S~4.2]{J}
(who states a loss of 1~bit $\approx$ 0.3 decimal place per step~$n$).
Still, the range (\ref{LIt}) needed for new tests of RH stays far out of reach.
All that motivates quests for simpler realizations of the Keiper--Li idea.
\smallskip

But first, as the \emph{asymptotic} sensitivity to RH is the main property we will prove to generalize,
we review that feature for $\{ \l _n \}$ itself.

\begin{figure}[h]
\center
\includegraphics[scale=.55]{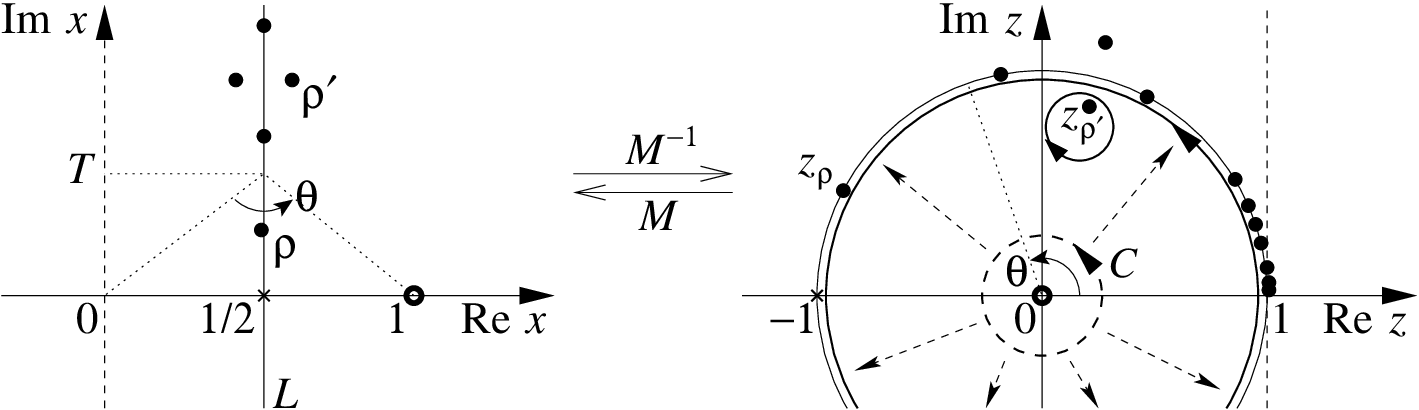}
\caption{\small Riemann zeros ({\footnotesize \textbullet}) depicted \emph{schematically} 
in the $x$ (left) and $z$ (right) upper half-planes
(at mock locations, including a putative pair off the critical line~$L$).
Symmetrical lower half-planes, and the zeros therein, are implied.
In the $z$-plane (right), we also plot the contour deformation used by Darboux's method 
upon the integral (\ref{Kinp}) for $n \to \infty$.
}
\end{figure}

\subsection{Asymptotic analysis of $\{ \l _n \}$ (to be generalized)}
\label{ADer}

The following derivation of the large-$n$ alternative (\ref{ASF})--(\ref{ASC}) for $\l _n$ 
readily settles the RH-false case (\ref{ASF}), 
then needs one more step for the RH-true formula (\ref{ASC}).
(In \cite{V} we obtained both cases in parallel by the saddle-point method used on a single integral,
but this approach does not yet extend.)

In this problem of large-order asymptotics, 
we can initially stick to Darboux's basic idea \cite[\S~7.2]{Di}:
for a sequence like (\ref{KEd}), its $n \to \infty$ form is ruled 
by the \emph{singularities} of its generating function $\varphi (z)$,
through the integral formula (equivalent to (\ref{KEd}) by the residue theorem)
\beq
\label{Kint}
\l ^{\rm K}_n = \frac{1}{2 \pi\mi} \oint_C \frac{\d z}{z^{n+1}} \varphi (z) , 
\qquad \varphi (z) \equiv \log 2 \xi \Bigl(\frac{1}{1-z}\Bigr) ,
\eeq
where $C$ is a positive contour close to $z=0$ leaving all other singularities 
(namely, those of $\varphi (z)$) outside.

\subsubsection{Darboux's method does the RH false case}
\label{DMF}

Solely for this stage, it is worth integrating (\ref{Kint}) by parts first, to
\beq
\label{Kinp}
\l ^{\rm L}_n = (2 \pi\mi)^{-1} \oint_C z^{-n} \, \frac{\d \varphi }{\d z} \, \d z :
\eeq
being meromorphic, $\d \varphi /\d z$ will be simpler to use than the multiply-valued~$\varphi $.

Since the integrand in (\ref{Kinp}) has the large-$n$ form $\e^{\Omega _n (z)}$ 
where $\Omega _n $ tends to~$\infty$ with $n$ ($\Omega _n (z) \sim -n \log z$), 
we may use the steepest-descent method \cite[\S~2.5]{Er}
to deform the contour $C$ toward decreasing $\Re (-n \log z)$, 
i.e., as a circle of radius $r$ growing toward~1 (Fig.~1 right);
then, by the residue theorem, each singularity of $\d \varphi /\d z$ swept in turn, 
namely a simple pole $z_{\rho '}$ per RH-violating zero $\rho '$ on the $\{ \Re x > \hf \}$ side,
yields an asymptotic contribution $(-z_{\rho '}^{\, -n})$ in descending order, 
and these altogether add up to (\ref{ASF}) \cite{V}.

If now RH is true, then as the radius of the contour attains $r=1^-$,
(\ref{ASF}) reaches no better than $\l ^{\rm L}_n = o(r^{-n}) \ (\forall r<1)$;
only a finer analysis of the limiting integral at $r=1^-$ 
pins down an asymptotic form for~$\l ^{\rm L}_n $, see next.

\subsubsection{Oesterl\'e's argument for the RH true case \cite{O}}
\label{KLO}
(reworded by us). Its starting point will be a real-integral form which comes from letting 
the contour in (\ref{Kint}) go up to $\{ |z|=1^- \}$ (unobstructed, under RH true),
making the change of variable $z = \e^{\mi \theta }$, 
and reducing to an integral over real $\theta $ with the help of the Functional Equation (\ref{CZ}):
namely, \cite{O}
\beq
\label{LR}
\l ^{\rm K}_n = 2 \int_0^\pi \sin n \theta \ N(\hf \cot \hf \theta ) \, \d \theta .
\eeq
Now, $\theta $ real $ \iff M(z) = \hf + \mi T$ with real $ T \equiv \hf \cot \hf \theta $,
so $\theta $ is also the angle subtended by the real vector $(\overrightarrow{01})$ 
from the point $\hf + \mi T$ (Fig.~1 left);
the counting function condenses to the critical line: $N(T) \equiv \# \{ \rho \in [\hf,\hf+\mi T] \}$.
An alternative validation of (\ref{LR}) is that its (Stieltjes) integral by parts
$n^{-1} \! \int_0^\pi 2(1-\cos n \theta) \, \d N$ at once yields the sum formula (\ref{KEp}) under~RH.
\smallskip

The $n \to \infty$ form mod $o(1)$ of (\ref{LR}) now directly stems 
from the Riemann--von Mangoldt large-$T$ form (\ref{RvM})
which, by virtue of $T \sim 1/\theta $, amounts to 
\beq
\label{NRv}
N(\hf \cot \hf \theta ) 
= - \frac{1}{2\pi \theta } \log (2\pi {\rm e} \, \theta) + O(|\log \theta |)_{\theta \to 0^+} ,
\eeq
that we plug into (\ref{LR}).

1) $O(|\log \theta |)_{\theta \to 0^+}$
is integrable in $\theta $ \emph{up to} $\theta =0 :$
by the Riemann--Lebesgue lemma \cite[\S~12.511]{GR} its integral against $\sin n \theta $ in (\ref{LR}) 
is $o(1)$, i.e., negligible;

2) as for the main term 
$- \int_0^\pi \sin n \theta \ (\pi \theta )^{-1} \log (2\pi {\rm e} \, \theta)\, \d \theta $ 
in (\ref{LR}), we change to the variable $\Theta _n \equiv n \theta $,
and then can push the new upper integration limit $n \pi$ to $+\infty$, 
mod $o(1)$ because the resulting integral is semiconvergent, to get
\[
\int_0^\infty \sin \Theta _n \, \frac{n}{\pi \Theta _n} 
(\log n - \log 2\pi {\rm e} \Theta _n ) \frac{\d \Theta _n}{n} .
\]
Now this is just $\ttfr{1}{\pi }[(\log n - \log 2\pi {\rm e})I_0-I_1]$, in terms of the classic integrals
\beq
\label{2CI}
I_0 = \int_0^\infty \sin \Theta \, \frac{\d \Theta }{\Theta } = \hf \pi , \qquad 
I_1 = \int_0^\infty \sin \Theta \, \log \Theta \, \frac{\d \Theta }{\Theta } = - \hf \pi \g 
\eeq
\cite[eqs.~(3.721(1)), (4.421(1))]{GR}. So that all in all,
\beq
\label{Kas}
\l ^{\rm K}_n = \hf \log n + c + o(1) , \quad c = \hf (\g -\log 2\pi -1) \qquad \mbox{(under RH true)},
\eeq
amounting to (\ref{OE}). \hfill $\square$

\section{An \emph{explicit} variant to the sequence $\{ \l ^{\rm K}_n \}$}

To alleviate the difficulties met with the original sequence $\{ \l _n \}$ (\S~\ref{PRH}),
we propose to \emph{deform} it (specifically, \emph{Keiper's} form (\ref{KEd})) 
into a simpler one, still RH-sensitive but of \emph{elementary closed form}.

While the original specification (\ref{KEd}) for $\l _n$ looks \emph{rigid}, 
this is only due to an extraneous assumption implied on the mapping $M$: 
that the origin $z=0$ \emph{has to} map to $x=1$ (the pole of $\z $).
Indeed, (\ref{KEd}) at once built $\{ \l _n \}$ upon the germ of $\log 2 \xi (x)$ at $x_0 \defi M(0)$ 
(the ``basepoint" for $\{ \l _n \}$), \emph{and} asserted $M(0) \equiv 1$ 
but \emph{this} was totally optional.
To wit, (\ref{KEd}) with a conformal mapping changed to $\widetilde M \ne M$
gives RH-sensitivity just as well provided (\ref{KLS}) holds,
and this only needs $\widetilde M ^{-1} (\{ \rho \} \cap L) \subset \{ |z|=1 \}$
which nowhere binds the basepoint $x_0 \ (= \widetilde M (0)$ now).
Prime examples are all $\widetilde M = M \circ H_{\tilde z}$, 
where $H_{\tilde z}$ conformally maps the unit circle onto itself as
\beq
\label{MT}
{z \mapsto H_{\tilde z}(z)} \defi {(z-\tilde z)/(1-\tilde z^\ast z)} \qquad 
(\tilde z \in {\mathbb C}, \ |\tilde z| \ne 1 );
\eeq
those $\widetilde M$, for which $x_0 = (1+\tilde z)^{-1}$, yield \emph{parametric} coefficients 
$\l _n (x_0)$ in terms of the derivatives $(\log \xi )^{(m)}(x_0)$;
for real $x_0$ excluding $\hf $,
these $\l _n (x_0)$ reproduce Sekatskii's ``generalized Li's sums" $k_{n,1-x_0}$ \cite{Se}.
Independently, different (double-valued) conformal mappings $\widetilde M$ generate ``centered" $\l ^0_n$ 
of basepoint precisely $x_0 = \hf $, the symmetry center for $\xi (x)$ (\cite[\S~3.4]{VK}, and Appendix).

But to attain truly simpler and explicit results, neither of those alterations goes far enough.
As a further step, rather than keeping a single basepoint (except, in a loose sense, $x_0=\infty \,$?),
we will crucially \emph{discretize} the derivatives of $\log \xi$ within the original $\l _n$ 
into selected \emph{finite differences} (and likewise in the Appendix for our centered $\l ^0_n$).

\subsection{Construction of a new sequence $\{ \L _n \}$}
\label{CNS}

The $\l _n$ are quite elusive as they involve \emph{derivatives} of $\log 2\xi $ 
(and worse, to any order), cf. (\ref{LId}).
On the integral form (\ref{Kint}), that clearly ties to the denominator $z^{n+1}$ 
having its zeros \emph{degenerate} (all at $z=0$, see fig.~2 left).

\begin{figure}[h]
\center
\includegraphics[scale=.55]{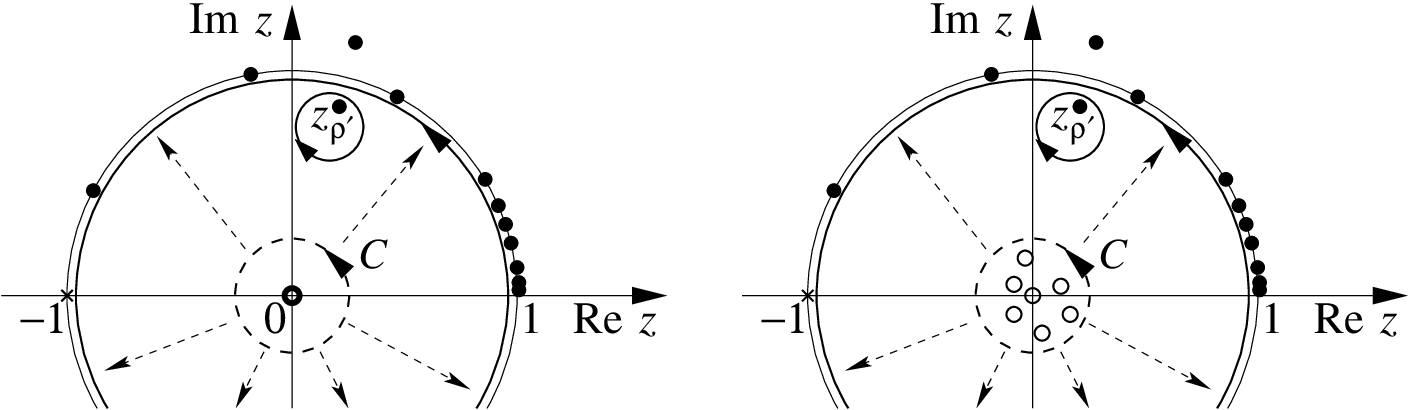}
\caption{\small As Fig.~1 right, but now showing the multiple pole $z=0$ (left plot) 
split into simple ones (right plot) without disrupting the large-$n$ analysis of \S~\ref{ADer}.
}
\end{figure}

Now at given $n$, if we \emph{split those zeros apart} as $0, z_1, \ldots, z_n$ 
(all distinct, and still inside the contour: fig.~2 right), then by the residue theorem, 
(\ref{Kint}) will simplify to a linear combination of the $\varphi (z_m)$, i.e., \emph{a finite difference}.
Doing so by plain shifts of the factors ($z \mapsto z-z_m$) would also split their unit disks apart,
and thus ruin the setting of \S~\ref{ADer} for asymptotic RH-sensitivity.
So to fix $\{ |z| \le 1 \}$, we use \emph{hyperbolic} shifts, i.e., 
again (\ref{MT}) (now with $|\tilde z|<1$), but \emph{differing on each factor}: $z \mapsto H_{z_m}(z)$ 
instead of one $H_{\tilde z}$ 
on all of $z^{n+1}$ previously (which left $\log 2 \xi$ differentiated all the same, only elsewhere).

The origin $z=0$ has then lost all special status, hence so does the particular mapping $M$ 
(selected to send $z=0$ to the pole $x=1$);
thereupon the variable~$x$, natural for the $\z $-function, is the simplest to use here as well. 
With everything rewritten in the variable $x$, (\ref{Kint}) reads as
\beq
\label{KLX}
\l ^{\rm K}_n = \frac{1}{2 \pi\mi} \oint \frac{\d x}{x(x-1)}
\Bigl( \frac{x}{x-1} \Bigl)^n \, \log 2\xi (x) \quad ( \mbox{integrated around }x=1) ,
\eeq
and the deformations just introduced have the form
\beq
\label{KLD}
\frac{1}{2 \pi\mi} \oint_{{\cal C}_n} \frac{\d x}{x(x-1)} 
\frac{1}{h_{x_1}(x) \ldots h_{x_n}(x)} \log 2\xi (x) , \ \quad 
h_{\tilde x} (x) \equiv \frac{\tilde x^\ast }{\tilde x} \, \frac{x-\tilde x}{x + \tilde x^\ast -1},
\eeq
($h_{\tilde x} \equiv H_{\tilde z} \circ M^{-1}, \ \tilde x \equiv M(\tilde z)$);
the contour ${\cal C}_n$ encircles the set $\{ 1,x_1,\ldots,x_n \}$ positively
(and may as well depend on~$n$). Then the integral in (\ref{KLD}) readily evaluates to
\beq
\label{Bint}
\sum_{m=1}^n \Res _{x=x_m} f_n (x) \log 2\xi (x_m) , \qquad 
f_n (x) \defi \frac{1}{x(x-1)} \, \frac{1}{[ h_{x_1} \ldots h_{x_n} ](x)} ,
\eeq
by the residue theorem ($x=1$ contributes zero thanks to $\log 2\xi (1) =0$).

\emph{Now, we pick ${x_m \equiv 2m}$ for $m=1,2,\ldots$} (independently of~$n$),
to later benefit from the known values $\z (2m)$. That fixes
\bea
\label{F_N}
f_n (x) \si=\se \frac{g_n(x)}{x(x-1)} \qquad (\sim 1/x^2 \quad \mbox{for }x \to \infty ) , \\
\label{G_N}
g_n(x) \si\defi\se \prod_{m=1}^n \frac{x + 2m - 1}{x-2m} \equiv 
\frac{ \G (\hf x-n) \, \G (\hf (x \!+\! 1) + n) }{ \G (\hf x) \, \G (\hf (x \!+\! 1)) } \\
\label{gdf}
\si\equiv\se g(x) (-1)^n \frac{\G (\hf (x \!+\! 1) + n)}{\G(1-\hf x + n)}, \qquad
g(x) \defi \frac{\sqrt \pi \, 2^{x-1}}{\sin (\pi x/2) \, \G (x)} 
\eea
(by the duplication and reflection formulae for $\G $; 
we refer to \S~\ref{MNB} for all current notations). The resulting residues of $f_n$ are
\[
\Res _{x=1} f_n = -\frac{(-1)^n }{A_{n0}}, \quad \Res _{x=2m} f_n = (-1)^{n+m} A_{nm} \quad 
(m=0,1,\ldots,n) ,
\]
\beq
\label{Anm}
A_{nm} = 2^{-2n} \frac{(n \!+\! m,n \!-\! m,2m)!}{2m-1} 
\equiv \frac{1}{(2m \!-\! 1) \, (n \!-\! m)! \, m!} \, \frac{\G (n \!+\! m \!+\! 1/2)}{\G (m \!+\! 1/2)} .
\eeq

Later we will need the partial fraction decomposition itself,
\beq
\label{PFD}
f_n (x) = (-1)^n \biggl[ \frac{-1}{A_{n0}} \, \frac{1}{x-1} + \sum_{m=0}^n \frac{(-1)^m A_{nm} }{x-2m} \biggr] :
\eeq
enforcing $f_n (x) \mathop{\sim}\limits _{x \to \infty} 1/x^2$ from (\ref{F_N}) 
first bars any additive constant but 0 in the right-hand side of (\ref{PFD}), as shown, then yields two more identities:
\beq
\label{AID}
\sum_{m=0}^n (-1)^m A_{nm} \equiv \frac{1}{A_{n0}} , \quad \mbox{resp.} \quad
2 \sum_{m=1}^n (-1)^m A_{nm} m \equiv (-1)^n + \frac{1}{A_{n0}} ,
\eeq
to be useful in \S~\ref{RKS} (these also exemplify sums of all the residues being $\equiv 0$
for meromorphic functions on the Riemann sphere: here $f_n(x)$, resp.~$xf_n(x)$).

Next, for each $n$ we select a positive contour ${\cal C}_n $ that wraps around 
the real interval $[1,2n]$ (the subinterval $[2,2n]$ would suffice, 
but here it will always prove beneficial to \emph{dilate}, not shrink, ${\cal C}_n $). 
Our final result is then
\bea
\label{Ldef}
\L _n \si\defi\se \frac{1}{2 \pi\mi} \oint_{{\cal C}_n } f_n(x) \log 2 \xi (x) \, \d x 
\qquad \mbox{(with $f_n$ from (\ref{F_N}))} \\
\label{EKL}
\si\equiv\se (-1)^n \sum_{m=1}^n (-1)^m A_{nm} \log 2 \xi (2m) , \qquad \qquad n=1,2,\ldots ,
\eea
and the latter form is \emph{fully explicit}: $A_{nm}$ are given by (\ref{Anm}), and 
\beq
\label{BE}
2 \xi (2m) = \frac{ (-1)^{m+1} B_{2m} }{ (2m-3)!! } \, (2\pi)^m \equiv 
\frac{ 2 (-1)^{m+1} B_ {2m} }{\G (m-\hf )} \, \pi^{m+1/2} , \quad m=0,1,\ldots 
\eeq
\beq
\label{ExV}
E.g., \ \L _1 = \frac{3}{2} \log \frac{\pi}{3}\, ,
\ \L _2 = \frac{5}{24} \log \biggl[ \Bigl( \frac{2}{5} \Bigr) ^{\! 7} \frac{3^{11}}{\pi^4} \biggr] ,
\ \L _3 = \frac{21}{80} \log \biggl[ \frac{5^{25} \pi ^8}{2^3 (3^2 \cdot 7)^{11}} \biggr] .
\eeq

So, we deformed Keiper's $\{ \l ^{\rm K}_n \}$ by discretizing the derivatives on $\log 2\xi $ 
to \emph{finite differences} anchored at locations $1, \{2m\}$ where $\xi $ has known values,
in a canonical way basically dictated by the \emph{preservation of RH-sensitivity}.\footnote
{One still has the freedom to spread the $x_m$ further out by skipping some locations~$(2j)$,
(or inversely, to keep some residual degeneracy in the numerics if ever that helps it).}
This gave $\{ \L _n \}$ the \emph{elementary closed form} (\ref{EKL}),
which is moreover directly computable at any individual~$n$ (in welcome contrast to the original~$\l _n$, 
which need an iterative procedure all the way up from $n=1$).

\subsection{Remarks.}
\label{RKS}

1) Thanks to the second sum rule (\ref{AID}),
the $(\log 2\pi)$-contributions to (\ref{EKL}) from the first expression (\ref{BE})
can be summed, resulting in $\L_n \equiv \hf \log 2\pi + u_n $ with
\beq
\label{Udef}
u_n \defi (-1)^n \left[ \sum_{m=1}^n (-1)^m A_{nm} \log \frac{ |B_{2m}|}{(2m-3)!!} 
+ \frac{1}{2 A_{n0}} \log 2\pi \right] :
\eeq
this sequence $\{ u_n \}$ is the one used in our first note \cite{VU}. Likewise, the rightmost expression (\ref{BE}) 
and the pair (\ref{AID}) lead to the partially summed form
\bea
\label{Vdef}
\L _n \si\equiv\se \hf \log \pi + (-1)^n 
\left[ \sum_{m=1}^n (-1)^m A_{nm} \log \frac{ |B_{2m}|}{\G (m \!-\! \hf )} \right. \nonumber \\
\si\se \qquad \qquad \qquad \quad \ \left. + \Bigl( \frac{1}{A_{n0}}\!-\! A_{n0} \Bigr) \log 2 
+ \Bigl( \frac{1}{A_{n0}}\!-\! \frac{A_{n0}}{2} \Bigr) \log \pi \vphantom{\sum_{m=1}^n} \right] ,
\eea
suitable for computer routines able to directly deliver $(\log \G )$.
\smallskip

2) If in place of (\ref{BE}) we evaluate 
\beq
\label{elr}
\log 2 \xi (2m) = \log [2(2m-1) m! \, \pi ^{-m}] + \sum_{p \ \rm prime} \, \sum_{r=1}^\infty \frac{p^{-2mr}}{r}
\eeq
(using (\ref{CZ}) and the expanded logarithm of the Euler product for $\z $),
then (\ref{EKL}) gives an \emph{arithmetic} form for $\L _n$, 
like Bombieri--Lagarias's for $\l ^{\rm L}_n$ \cite[Thm~2]{BL}. 
For numerics moreover, (\ref{elr}) becomes \emph{increasingly efficient as $m$~grows}, 
exactly counter to (\ref{BE}). And the corresponding partially summed form~is
\bea
\L _n \si\equiv\se -\hf \log \pi + (-1)^n 
\! \left[ \sum_{m=1}^n (-1)^m A_{nm} \biggl[ \log ((2m \!-\!1)m!) 
+ \! \sum_{p \ \rm prime} \, \sum_{r=1}^\infty \frac{p^{-2mr}}{r} \biggr] \right. \nonumber \\[-4pt]
\si\se \qquad \qquad \qquad \qquad \qquad \qquad \quad \left. 
+ \Bigl( \frac{1}{A_{n0}}\!-\! A_{n0} \Bigr) \log 2 
- \frac{1}{2A_{n0}} \log \pi \vphantom{\sum_{m=1}^n} \right] .
\eea
\smallskip

3) B\'aez-Duarte \cite{BD} has an equally explicit sequential criterion for RH in terms of $\{ B_{2m} \}$, 
but its detection domain for RH-violations above $T_0$ is \emph{exponentially} distant, 
$n \gtrsim \e ^{\pi T_0}$ \cite[\S~4]{M3} (\cite[\S~7]{FV} quotes $n \gtrsim 10^{600,000,000}$);
for our $\L _n$ the corresponding $n$-threshold will prove considerably closer (\S~\ref{FRH}).

4) The whole scheme will be extended from the Riemann zeta function~$\z $
to certain \emph{Dirichlet L-functions} in \S~\ref{DLF}, then to some linear combinations thereof, 
specifically the \emph{Davenport--Heilbronn} functions, in \S~\ref{DHF}.

\subsection{Expression of $\L _n$ in terms of the Riemann zeros}
\label{LRZ}

Let the primitive $F_n$ of the function $f_n$ in (\ref{F_N}), (\ref{PFD}) be defined by
\bea
\label{Fdef}
F_n(x) \si\defi\se \int_\infty^x f_n(y)\, \d y \qquad \qquad 
(\Rightarrow F_n (x) \sim -1/x \mbox{ for } x \to \infty ) \nonumber \\
\si\equiv\se (-1)^n \biggl[ - \frac{1}{A_{n0}} \log (x-1)
+ \sum_{m=0}^n (-1)^m A_{nm} \log (x-2m) \biggr] 
\eea
and by single-valuedness in the whole $x$-plane minus the cut $[0,2n] :$
e.g., $F_1(x) = {\hf \log \, [x(x-2)^3/(x-1)^4]}$.

Then in terms of (\ref{Fdef}), the $\L _n$ are expressible as sums over the Riemann zeros 
(which converge like $\sum_\rho 1/\rho $ for any $n$, hence need the rule (\ref{SG})):
\beq
\label{LIR}
\L _n \equiv \sum_\rho F_n (\rho) , \qquad n=1,2,\ldots .
\eeq
(In the original $\l ^{\rm K}_n $, (\ref{KLX}) uses $[x/(x-1)]^n$ in place of $g_n (x)$, 
which exceptionally yields \emph{rational} functions $n^{-1}[1-(1-1/(1-x))^n]$ in place of the $F_n(x)$, 
for which (\ref{LIR}) restores~(\ref{KEp}).)

\begin{figure}[h]
\center
\includegraphics[scale=.6]{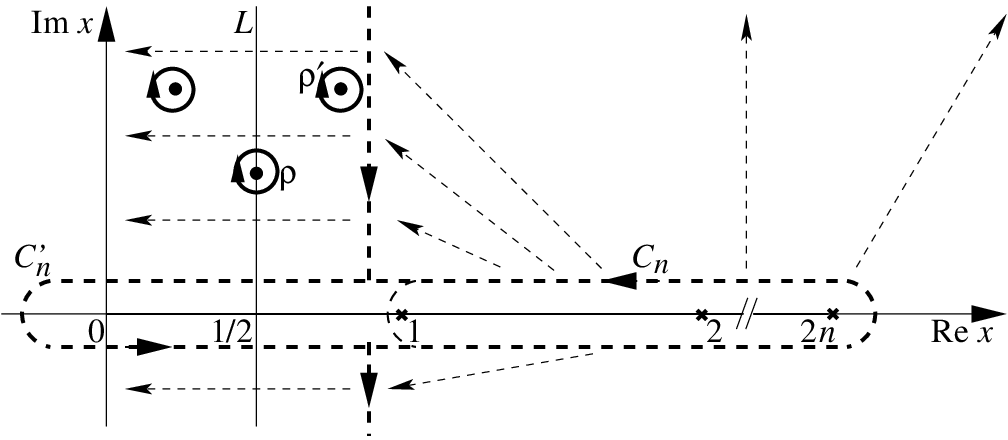}
\caption{\small As Fig. 1 left, but superposing the deformation of the integration path for the integral (\ref{Lfed}) 
in which the function $\ttfr{1}{\xi } \ttfr{\d \xi }{\d x}$ has the Riemann zeros as simple poles.
A symmetrical lower half-plane is implied.}
\end{figure}

Proof of (\ref{LIR}) (outlined, see fig. 3): first stretch the contour ${\cal C}_n$ in (\ref{Ldef}) 
to ${\cal C}_n '$ fully enclosing the cut $[0,2n]$ of $F_n$ (as allowed by $\log 2\xi (0)=0$). 
As $F_n$ is single-valued on ${\mathbb C} \setminus [0,2n]$, 
the so modified (\ref{Ldef}) can be integrated by parts,
\beq
\label{Lfed}
\L _n \defi -\frac{1}{2 \pi\mi} \oint_{{\mathcal C}_n '} F_n(x) \frac{1}{\xi } \frac{\d \xi }{\d x} \, \d x ,
\eeq
then the contour ${\mathcal C}_n '$ can be further deformed into a sum of 
an outer anticlockwise circle ${\mathcal C}_R$ centered at $\hf $ of radius $R \to \infty$ (not drawn),
and of small clockwise circles around the poles of the meromorphic function 
$\ttfr{1}{\xi } \ttfr{\d \xi }{\d x}$ inside~${\mathcal C}_R$; 
these poles are the Riemann zeros $\rho $ therein, and each contributes $F_n(\rho )$.
By the Functional Equation (\ref{CZ}), the integral on ${\mathcal C}_R$ is also
$\oint_{{\mathcal C}_R} \langle F_n (x) \rangle \, \ttfr{1}{\xi } \ttfr{\d \xi }{\d x} \, \d x$, 
where $\langle F_n (x) \rangle \defi {\hf [F_n (x)+F_n (1 \!-\! x)]} =O(1/x^2)$;
hence this integral tends to~0 if $R \to \infty$ while keeping far enough
from ordinates of Riemann zeros in a classic fashion 
(i.e., so that $\Bigl| \ttfr{1}{\z } \ttfr{\d \z }{\d x} \Bigr| (r + \mi R) < K \log ^2 R$ 
for all $r \in [-1,+2]$, cf. \cite[p.~108]{Da}), hence (\ref{LIR}) results. \hfill $\square$

\section{Resulting new sequential criterion for RH}

Like the elusive Keiper--Li sequence, the fully explicit one $\{ \L _n \}$ proves RH-sensitive 
(just slightly differently). The argument will use the function
\bea
\label{MOM}
\Phi (X) \si\defi\se X \log \Bigl( 1- \frac{1}{X^2} \Bigr) + \log \frac{X+1}{X-1}
\equiv \int _\infty ^X \log \Bigl( 1- \frac{1}{Y^2} \Bigr) \d Y \\ 
\label{PD}
\si\equiv\se (X+1) \log(X+1) + (X-1) \log(X-1) -2 \, X \log X 
\eea
made single-valued on ${\mathbb C} \setminus [-1,1] \cup \{ \infty \}$
and \emph{odd}, so that $\mi ^{-1} \Phi $ becomes \emph{real} on~$\mi {\mathbb R}$:
\beq
\label{PsD}
\Psi (U) \defi \mi ^{-1} \Phi (\mi U) \equiv U \log (1+1/U^2) -2 \arctan(1/U) , \quad 
U \in {\mathbb R}^\ast \cup \{ \infty \}.
\eeq 

\subsection{Asymptotic criterion}

Our main result is an \emph{asymptotic} sensitivity to RH as $n \to \infty$,
through this alternative for $\{ \L _n \}$ which parallels (\ref{ASF})--(\ref{ASC}) for $\{ \l _n \} :$
\bea
\label{LNRH}
\bullet \ \mbox{RH false:} \ \L _n \si\sim\se 
\sum _{\rho ' - \frac{1}{2} \, \in \, 2n \, {\mathcal D}_{R_0} } F_n (\rho ') \quad \pmod{o(R^n) \ \forall R > R_0 >1 } , \\
\mbox{where} \quad {\mathcal D}_{R_0} \si\defi\se \{ X \in {\mathbb C} \mid \Re \Phi (X) > \log R_0 \} \nonumber\\
\si\subset\se \{ \Re X >0 \},
\mbox{ so the sum (\ref{LNRH}) is a truncation of} \sum _{\Re \rho ' > \frac{1}{2} } , \nonumber\\[-8pt]
\label{ST1}
\mbox{and} \quad F_n (\rho ') \si\sim\se 
\frac{g(\rho ')}{\rho '(\rho '-1)} (-1)^n \frac{n^{\rho '-\frac{1}{2} }}{\log n} \quad 
(n \to \infty ) \quad \mbox{for each given } \rho ' \\
\si\se \mbox{(leads to a growing oscillation \emph{with both signs}, cf. (\ref{gdf}));} \nonumber\\
\label{LRH}
\bullet \ \mbox{RH true:} \ \L _n \si\sim\se 
\log n + C \pmod{o(1)} , \qquad C = \hf(\g - \log \pi - 1) \\
\si\se \mbox{(implying asymptotic \emph{positivity})} \nonumber
\eea
which, compared with (\ref{Kas}), yields $C \equiv c+\hf \log 2 \ (\approx -0.783757110474)$, 
and testifies that we kept $\L _n$ qualitatively close to $\l ^{\rm K}_n$.

For the RH false case, the pair (\ref{LNRH})--(\ref{ST1}) is formally, in the variable $\log n$, 
an expansion in exponentials multiplied by divergent power series (i.e., a \emph{transseries}), 
to be interpreted with caveats (detailed in \S~\ref{RHF} below).
\medskip

The derivation scheme will transpose the arguments of \S~\ref{ADer} to $\L _n$ expressed in an integral form 
and now in the $x$-plane: by either (\ref{Ldef}), using the function $g_n$, in place of (\ref{Kint});
or its integral by parts (\ref{Lfed}), using $F_n$, in place of (\ref{Kinp}). 
Two new problems arise here: the large-$n$ forms of $g_n (x)$ and $F_n (x)$ need to be worked out,
and the geometry of the integrands is strongly $n$-dependent hence the relative scales
of $n$ and $x$ will matter.

\textbf{3.1.1 If we need \emph{uniform} asymptotics in the integrands:} we must rescale the geometry as, for instance, 
$x = \hf + 2nX :$ this condenses the singularities onto the fixed $X$-segment $[0,1]$ as $n \to \infty$.
The Stirling formula applied to the ratio of $\G$ functions in (\ref{G_N}) then brings in
the function (\ref{MOM}):
\beq
\label{gAS}
g_n (\hf + 2nX) \sim \biggl( \frac{X+1}{X-1} \biggr) ^{1/4} \e^{n \Phi (X)} \qquad \mbox{for } n \to \infty .
\eeq

Next, integrating for $F_n (x) = \displaystyle \int _\infty ^x \frac{g_n (y)}{y(y-1)} \d y$, 
with $y = \hf + 2nY$ yields
\beq
F_n (\hf + 2nX) \sim \frac{1}{2n} 
\int _\infty ^X \frac{\d Y}{Y^2} \biggl( \frac{Y+1}{Y-1} \biggr) ^{1/4} \e ^{n \Phi(Y)},
\eeq
which is a Laplace transform in the integration variable $\Phi $, hence it has an asymptotic power series in $(1/n)$ 
(usually divergent) starting as \cite[eq.~2.2(2)]{Er}
\bea
F_n (\hf + 2nX) \si\sim\se
\frac{1}{2n} \, \frac{1}{X^2} \biggl( \frac{X+1}{X-1} \biggr) ^{1/4} \frac{ \e ^{n \Phi(X)}}{n \ \d\Phi /\d X} \nonumber\\
\label{FAN}
\si=\se 
\frac{1}{2n^2} \, \frac{1}{X^2} \biggl( \frac{X+1}{X-1} \biggr) ^{1/4} \frac{ \e ^{n \Phi(X)}}{\log (1-1/X^2)}
\eea

\textbf{3.1.2 Whereas if we let $n \to \infty$ \emph{at fixed} $x$:} (\ref{gdf}) at once implies
\beq
\label{gAs}
g_n (x) \sim g(x) (-1)^n n^{x-1/2} \sim g(x) (-1)^n \e^{\log n \, (x-1/2)} .
\eeq
Here $(\log n)$ replaces $n$ as large parameter, and the same Laplace argument for the integral 
$F_n (x) = \displaystyle \int _\infty ^x \frac{g_n (y)}{y(y-1)} \d y$
now yields an asymptotic series in powers of $(1/\log n)$ (usually divergent), starting as
\beq
\label{FAS}
F_n (x) \sim \frac{g(x)}{x(x-1)} (-1)^n \frac{n^{x-1/2}}{\log n} .
\eeq

\subsection{Details for the case RH false}
\label{RHF}

As in \S~\ref{DMF}, we can apply the steepest-descent method to the integral (\ref{Lfed}) written in the global variable $X$.
We rescale the contour ${\mathcal C}_n '$ then deform it toward level contours $\Re \Phi (X) = \Phi _0 \to 0^+$
which thus approach the completed critical line in the $X$-plane, $\{ \Re X =0 \} \cup \{ \infty \}$, 
from the ${\{ \Re X >0 \}}$ side (fig.~4). 
Apart from staying on this side (and being rescaled to the $X$-plane), 
the contour deformation is isotopic to that used in \S~\ref{LRZ},
hence it likewise yields a contribution $F_n (\rho ')$ per RH-violating zero $\rho '$ whose 
rescaled image $X_{\rho '} \equiv (\rho ' -\hf ) / (2n)$ has $\Re \Phi (X_{\rho '}) > \Phi _0 >0$;
i.e., overall,
\[
\sum _{\Re \Phi (X_{\rho '}) > \Phi _0} F_n (\rho ') .
\]
(The novelty vs (\ref{LIR}) is that now the terms are ordered asymptotically.)
For ${X \to \infty}$, $\Phi (X) \sim 1/X$, hence the level curves asymptotically become $\{ \Re (1/X)= \Phi _0 \}$ 
(circles tangent to the imaginary $X$-axis at 0).
Thus the above sum over zeros has a natural cutoff $|\Im \rho '| \lesssim \sqrt{n/ \Phi _0}$, 
which is the height above which the disk $\{ \Re (1/X) > \Phi _0 \}$ parts from the critical strip.

\begin{figure}[h]
\includegraphics[scale=.6]{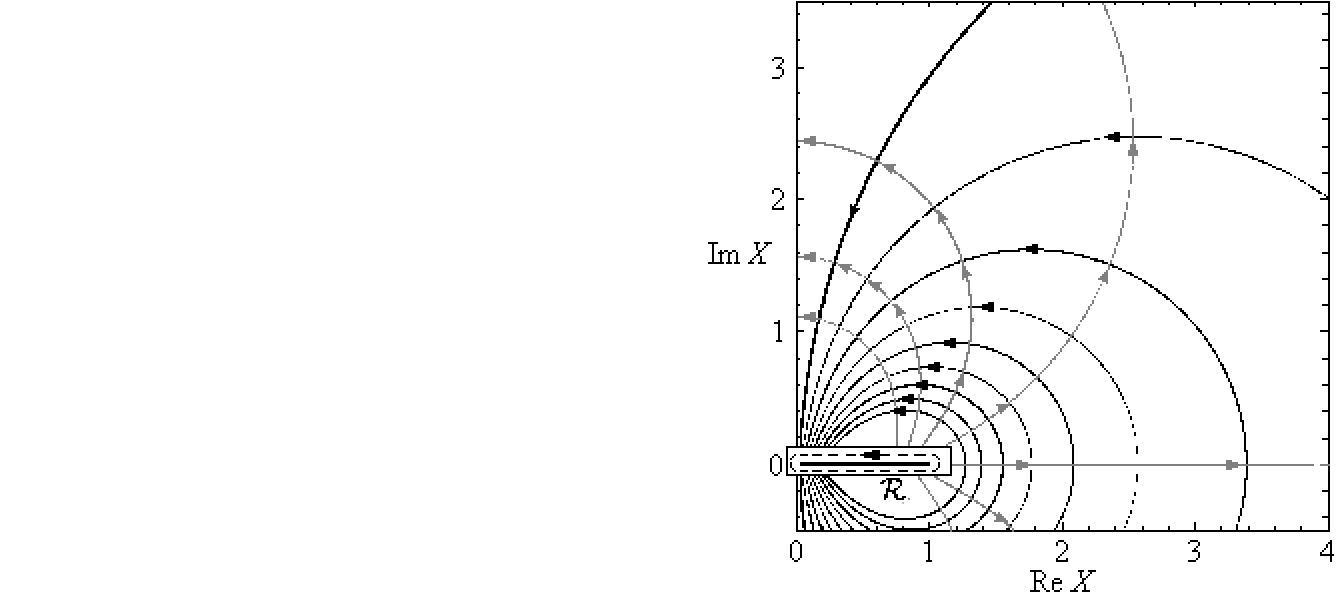}
\vskip -3.8cm
\includegraphics[scale=.4]{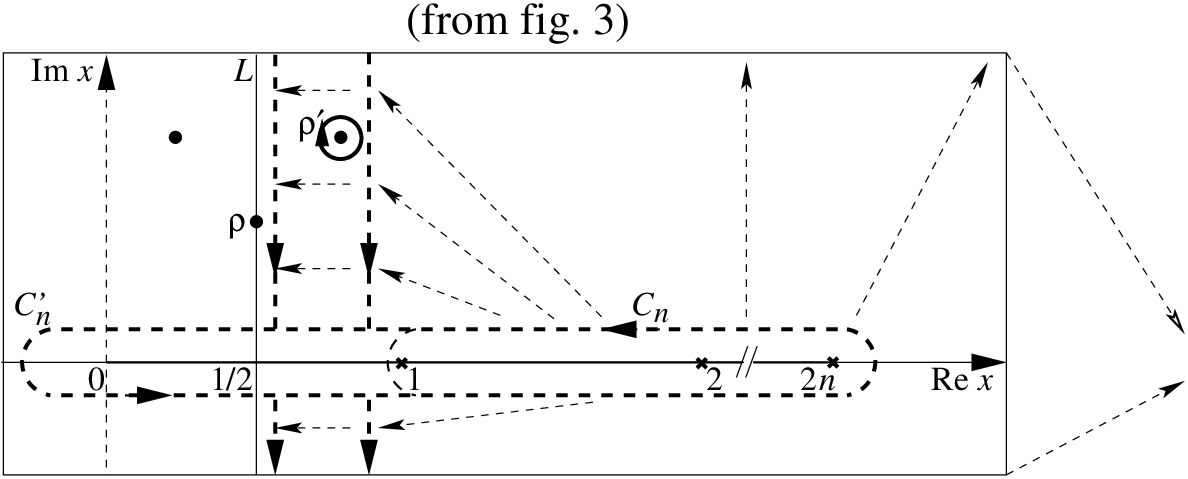}
\vskip 4mm
\caption{\small (at right) as fig. 3 but using the global variable $X \defi (x-\hf )/(2n)$ 
in the integral (\ref{Lfed}) (fig.~3 now fills the rectangle $\mathcal R$). 
Mathematica \protect\cite{W} contour plots for $\Phi (X)$ of (\ref{MOM}): 
in black, deformed integration contours (level curves ${ \{ \Re \Phi (X) = \Phi _0 \} }$, 
$\Phi _0 \downarrow 0^+$); in gray, steepest-descent lines (level curves of $\Im \Phi $).}
\end{figure}

When zooming in to the $n \to \infty$, fixed-$x$ regime, then by (\ref{FAS}) the level curves turn to 
parallel lines $\{ \Re x - \hf = t_0 \}$, and the deformation becomes ${t_0 \to 0^+}$ (fig.~4 left).
Only a portion of contour ${\mathcal C}_n ' \cap {\{ \Re x -\hf <\eps \} }$ escapes the deformation, 
but by (\ref{FAS}) its contribution to the integral (\ref{Lfed}) is $O(n^\eps )$, ultimately negligible.
Therefore, $\L _n \pmod{o(R^n) \ \forall R>\e^{\Phi _0}}$ equals the above sum $\sum_{\rho '} F_n(\rho ')$, 
which proves (\ref{LNRH}) under the replacement $R_0 = \e^{\Phi _0}$.

This fixed-$x$ regime actually governs \emph{each individual term} of the sum~(\ref{LNRH}):
$F_n (\rho ')$  is the contribution to $\L _n $ by a given zero $\rho '$, that lives at a fixed~$x$
($\equiv \hf +t+\mi T \ (t>0)$).
Therefore, setting $x = \rho '$ in (\ref{FAS}) directly delivers the asymptotic form (\ref{ST1}) for $F_n (\rho ')$  
which is \emph{explicit} including the ${n\to \infty }$ dependence, 
as opposed to the exact $F_n (\rho ')$ specified via (\ref{Fdef}).
Now (\ref{LNRH}), which is our other asymptotic ($n \to \infty$) statement, still uses \emph{exact} $F_n (\rho ')$
that cannot be substituted outright by their explicit forms (\ref{ST1}), 
because (\ref{ST1}) is not asymptotic \emph{uniformly in} $T$; 
it is nevertheless what describes $F_n (\rho ')$ for $n \to \infty$. 
Even more explicitly, using $|T| > T_0 \gg 1$, 
$g(\rho ') \approx (2/|T|)^t \e^{-{\rm sign} (T) \mi\pi/4} (2 {\rm e}/|T|)^{\mi T}$ 
(from the Stirling formula) leads to
\beq
\label{LABS}
F_n(\rho ') \approx \frac{1}{T^2 \log n} \Bigl( \frac{2n}{|T|} \Bigr) ^t 
(-1)^{n+1} \e^{-{\rm sign} (T) \, \mi \pi /4} \Bigl( \frac{2 {\rm e} n}{|T|} \Bigr) ^{\mi T}
\quad \mbox{for } n \gg |T| .
\eeq

Thus, all in all, we \emph{needed} two separate asymptotic formulae 
((\ref{LNRH})--(\ref{ST1})) to fully describe $\{ \L _n \}$ in the RH-false case. 
\smallskip

If now RH is true then, as $\Phi _0$ attains $0^+$, (\ref{LNRH}) reaches no better than $\L _n = o(R^n) \ (\forall R>1)$,
and only a finer analysis of the limiting integral on the critical line $L$
will lead to a definite asymptotic form, as follows.

\subsection{Details for the case RH true}
\label{PTL}

Here our quickest path is to adapt Oesterl\'e's argument from \S \ref{KLO}.

To deform $\{ \l _n \}$ into $\{ \L _n \}$, we replaced $[ x/(x-1) ]^n$ 
in (\ref{KLX}) by $g_n(x)$ from~(\ref{G_N}). That changes (\ref{LR}) to
\beq
\label{GR}
\L _n = 2 \int_0^\pi \sin \Theta _n (\theta) \ N(T(\theta )) \, \d \theta , 
\qquad T(\theta ) \equiv \hf \cot (\hf \theta ) ,
\eeq
where $\Theta _n (\theta) \in (0,n\pi]$ ($\equiv n \theta $ previously for $\l _n$) is now the sum 
of the $n$~angles subtended by the real vectors $(\overrightarrow{1-2m,2m})$ from the point $\hf+\mi T (\theta)$, 
for $m=1,2,\ldots,n$. Namely,
\beq
\label{TDef}
\Theta _n (\theta) =2\sum_{m=1}^n \arctan \Bigl[ (4m-1) \tan (\hf \theta ) \Bigr] : 
\quad [0,\pi] \mapsto [0,n\pi] .
\eeq
The two \emph{endpoint slopes} of the function $\Theta _n (\theta )$ will then play key roles:
\bea
\label{TP0}
\Theta _n'(0) \si=\se \sum_{m=1}^n (4m-1) \equiv n(2n+1) , \\
\label{TPp}
\Theta _n'(\pi) \si=\se 
\sum_{m=1}^n (4m-1)^{-1} \equiv \qt \bigl[ \psi (n+\tq ) + \g + 3 \log 2 - \pi/2 \bigr] \\
\label{TPm}
\si=\se \min_{[0,\pi ]} \Theta _n'(\theta) \qquad \mbox{ because } 
\Theta _n '(\theta) >0, \ \Theta _n ''(\theta) <0 \mbox{ on } [0,\pi ] 
\eea
(all resulting from (\ref{TDef}) term by term;
we denote $' \equiv \ttfr{\d }{\d \theta }$, $\psi (x) \equiv \ttfr{\d }{\d x}\log \G (x)$).

We then follow the same steps as with (\ref{LR}) for $\l ^{\rm K}_n$ in \S~\ref{KLO}, 
now plugging $N(\hf \cot \hf \theta ) = -(2\pi \theta )^{-1} \log (2\pi {\rm e} \, \theta) 
+ O(|\log \theta |)_{\theta \to 0^+}$ from (\ref{NRv}) into (\ref{GR}).

1) $O(|\log \theta |)$ will be negligible as before if the derivation of the Riemann--Lebesgue lemma \cite[\S~5.14]{Ru}
extends so as to prove
\beq
\label{RLl}
I_n(f) \defi \int_0^\pi \sin \Theta _n(\theta ) \, f(\theta ) \, \d \theta = o(1)_{n \to \infty} 
\ \mbox{ for any $f$ integrable on }[0,\pi ] . 
\eeq 
Now $C^1$ functions are dense in $L^1$, \cite[thm 3.14]{Ru}
so for any $\eps >0$ there exists ${\mathbf f} \in C^1[0,\pi]$ such that $\Vert f-{\mathbf f} \Vert_{L^1} < \eps $, 
then it remains to prove $|I_n({\mathbf f})|<\eps $ for $n$ large enough, i.e., $I_n({\mathbf f})=o(1)_{n \to \infty}$. 
We integrate by parts in the variable~$\Theta _n$:
\bea
\label{IP}
I_n({\mathbf f}) \si=\se 
\int_0^{n\pi} \!\! \sin \Theta _n \Bigr[ {\mathbf f}(\theta (\Theta _n)) \frac{\d \theta }{\d \Theta _n} \Bigr] \d \Theta _n 
= \Bigr[ \!\! -\cos \Theta _n \, {\mathbf f}(\theta (\Theta _n)) \frac{\d \theta }{\d \Theta _n} \Bigr] 
_0^{n\pi } \! + J_n , \\
J_n \si\defi\se 
\int_0^{n\pi} \!\! \cos \Theta _n \, \frac{\d }{\d \Theta _n} 
\Bigr[ {\mathbf f}(\theta (\Theta _n)) \frac{\d \theta }{\d \Theta _n} \Bigr] \d \Theta _n
= J_n^{(1)} + J_n^{(2)} , \nonumber \\
J_n^{(1)} \si=\se 
\int_0^{n\pi} \!\! \cos \Theta _n \, {\mathbf f}'(\theta (\Theta _n )) 
\Bigl( \frac{\d \theta }{\d \Theta _n} \Bigr) ^{\! 2} \d \Theta _n \equiv 
\int_0^\pi \! \cos \Theta _n (\theta ) \, {\mathbf f}'(\theta ) \frac{1}{\Theta '_n(\theta )} \, \d \theta ,
\nonumber \\
J_n^{(2)} \si=\se 
\int_0^{n\pi} \!\! \cos \Theta _n \, {\mathbf f}(\theta (\Theta _n )) \frac{\d^2 \theta }{\d \Theta_n^2} \, \d \Theta _n
= \cos \Theta _n^{(0)} \, {\mathbf f}(\theta (\Theta _n^{(0)} )) \Bigl[ \frac{\d \theta }{\d \Theta _n} \Bigr]
_0^{n\pi } \! , \ \Theta _n^{(0)} \in [0,\pi] \nonumber
\eea
for \emph{some} $\Theta _n^{(0)}$, by the first mean-value theorem \cite[\S~12.111]{GR} (with 
$\d ^2 \theta / \d \Theta _n^2 \equiv -\Theta _n''/{\Theta _n'}^3$ \emph{not changing sign}, from (\ref{TPm})).
Now (\ref{TPp}) entails $\Theta _n'(\pi ) \sim \qt \log n$ upon which (\ref{TPm}) implies 
$\d \theta / \d \Theta _n \equiv 1/\Theta _n' (\theta) = O(1/\log n)$ uniformly.
Then $J_n^{(1)}, J_n^{((2)}$ are $O(1/\log n)$ from their respective rightmost expressions above,
hence so are $J_n$ and finally $I_n({\mathbf f})$ in (\ref{IP}), which ensures $I_n({\mathbf f})=o(1)$.
\hfill $\square$
\smallskip

2) For the main term 
$- \int_0^\pi \sin \Theta _n (\theta ) \, (\pi \theta )^{-1} \log (2\pi {\rm e} \, \theta) \, \d \theta $ in (\ref{GR}): 
by (\ref{RLl}), only the nonintegrable singularity multiplying $\sin \Theta _n $ 
should matter (i.e., ${\theta =0}$);
here $\Theta _n (\theta ) \sim \Theta _n'(0) \, \theta$, vs $n \theta$ previously, 
so it must suffice to substitute $\Theta _n'(0)=n(2n+1)$ for~$n$ 
in the previous result (\ref{Kas}) for $\l _n^{\rm K}$, to get
\beq
\label{LNA}
\L _n \sim \hf \log \Theta _n'(0) + c = \hf \log \, [n(2n+1)] + c \sim \log n + (c + \hf \log 2) .
\eeq
And the intermediate results to replace (\ref{2CI}) must be
\bea
\label{TH_N}
{\mathcal I}_{0,n} \si=\se \int_0^\pi \sin \Theta _n \frac{\d \theta }{\theta } 
= \hf \pi + o(1)_{n \to \infty} \, , \\
\label{TL_N}
{\mathcal I}_{1,n} \si=\se \int_0^\pi \sin \Theta _n \log \theta \, \frac{\d \theta }{\theta }  
= -\hf \pi \, (\log \, [n(2n+1)] + \g ) + o(1)_{n \to \infty} \, ,
\eea
which can be checked numerically (through nonarithmetical tests),
and indeed combine in (\ref{GR}) to yield (\ref{LNA}): i.e., the targeted result (\ref{LRH}).

Now a technicality arises: (\ref{G_N}) obviously implies 
\beq
\label{Thg}
\e ^{- \mi \Theta _n (\theta )} \equiv g_n (\hf +\mi T(\theta )) 
\qquad (\mbox{recalling } T(\theta ) \equiv \hf \cot \hf \theta ) ,
\eeq
so the function $\Theta _n $ also has the two disjoint nonuniform $n \to \infty$ regimes of~$g_n$: 
(\ref{gAS}) for $\theta = O(1/n)$ vs (\ref{gAs}) for fixed but \emph{nonzero} $\theta $,
and this complicates a rigorous rewording of the above heuristics, as done next. 
(We have quicker proofs for (\ref{TH_N}) alone.) 

We then shift focus onto a parametric integral, 
${\mathcal I}_n(\alpha ) \defi \int_0^\pi \sin \Theta _n \, \d \theta /\theta ^\alpha $ (for $\alpha <1$),
which in the end will restore ${\mathcal I}_{q,n} \equiv (-\d /\d \alpha )^q \, {\mathcal I}_n (\alpha =1^-)$ for $q=0,1$. 
But to better benefit from complex contour deformations, we prefer integrals \emph{without endpoints}:
with $\theta $~being an angle, such are 
\beq
\label{Ipm}
{\mathcal I}_{\pm ,n}(\alpha ) \defi 
\int_{-\pi }^\pi g_n (\hf +\mi T(\theta )) (\e ^{\pm \mi \pi /2} \theta )^{-\alpha } \, \d\theta ,
\eeq
where the standard determination (cut along ${\mathbb R}^-$) is taken for $v \mapsto v^{-\alpha }$.
The identity (\ref{Thg}) extends to $[-\pi ,\pi]$ with $\Theta _n (-\theta) = -\Theta _n (\theta)$,
which  actually makes $\e ^{- \mi \Theta _n (\cdot )}$ a \emph{regular} function 
of $\theta \in {\mathbb R} /2\pi {\mathbb Z}$, hence also $g_n (\hf +\mi T(\cdot ))$. 
Then by inspection, ${\mathcal I}_{\pm ,n}(\alpha ) \equiv 
2 \int _0^\pi \cos [\Theta _n (\theta) \pm \hf \pi \alpha ] \,  \d \theta /\theta ^\alpha $,
in turn implying
\beq
\label{IDD}
{\mathcal I}_n (\alpha ) \equiv 
\frac{1}{4 \sin \hf \pi \alpha } [{\mathcal I}_{-,n}(\alpha )  - {\mathcal I}_{+,n}(\alpha ) ] .
\eeq
Now (\ref{Ipm}) is controllable for $n \to \infty$ by the steepest-descent method: 
we deform its real integration cycle ($\{ |g_n|=1 \}$) to level lines 
${\{ |g_n| = {\rm const.} <1 \}}$ in the complex domain,
$\{ \theta - \mi u_\theta \mid \theta \in {\mathbb R} /2\pi {\mathbb Z} \}$,
$u_\theta >0$ due to $\Theta '_n (\theta) >0$.
(This is, viewed in the $\theta$-variable, just the deformation of \S~\ref{RHF} 
but taken beyond the imaginary $X$-axis.)
Since moreover  $\Theta '_n (\theta) \to +\infty$ uniformly, \emph{barely deformed} level curves 
(such that $u_\theta \sim \sqrt{\log n}/\Theta '_n (\theta )=o(1)$ uniformly) reach $|g_n|=o(1)$ 
save for obstructions to that lowering of the path, of which there may be three.
First, singularities of the function $\Theta _n (\theta )$, 
but the nearest one is $(-2 \arctan [1/(4n-1)]\ \mi )$, and $|g_n|=o(1)$ is achieved far above. 
Then, the jump we call $h(u)$ of the nonperiodic factor $\theta ^{-\alpha }$ in (\ref{Ipm}) 
across the line ${\{ \pi -\mi u \pmod {2\pi} \mid u \in {\mathbb R} \}}$:
its contribution to ${\mathcal I}_{\pm ,n}(\alpha )$ has size 
$| \! \int_0^{u_\pi } g_n (\hf +\mi T(\pi - \mi u)) \, h(u) \, \d u \, | \sim 
| \! \int_0^{u_\pi } g(\hf - \qt u) \, h(u) \, n^{-u/4} \d u \, |$
using ${T(\pi - \mi u)} \equiv \hf \cot \hf (\pi -\mi u) \sim \qt \mi u$ then (\ref{gAs})
(the regime of \S~3.1.2 holds as $T \approx 0$);
now both $g(\cdot)$ and $h$ are \emph{regular} functions in $u$ hence the last integral, 
by the \emph{Laplace method as in} \S~3.1.2, is $O(1/\log n) =o(1)$. 
Finally, the sole obstacle of consequence is the cut of $(\pm \mi \theta )^{-\alpha }$ from $\theta =0$ 
but for $(+ \mi \theta )^{-\alpha }$, even that cut goes \emph{upwards} (out of the way), 
resulting in ${\mathcal I}_{+,n}(\alpha )=o(1)$.
Whereas for ${\mathcal I}_{-,n}(\alpha )$, the cut is \emph{downward}: 
${ \{ \theta = -\mi u \mid u>0 \} }$, and it contributes
\beq
\label{OUF}
{\mathcal I}_{-,n}(\alpha ) \sim (\e ^{+\mi \pi \alpha } - \e ^{-\mi \pi \alpha }) 
\int_0^{u_0} \e ^{-\mi \Theta _n (-\mi u)} \frac{(-\mi \, \d u)}{u^\alpha }
= 2 \sin \pi \alpha \int_0^{u_0} \e ^{-\mi \Theta _n (-\mi u)} \frac{\d u}{u^\alpha } .
\eeq
Since $|u| \le u_0 \sim \ttfr{\sqrt{\log n}}{n(2n+1)}=o(1/n)$, we rescale $\tau \equiv 2nu \ (=o(1))$.
Now from~(\ref{Thg}), $ \e ^{-\mi \Theta _n (-\mi u)} 
\sim g_n(\hf - 2n/\tau ) \sim (\frac{1-\tau }{1+\tau })^{1/4} \e ^{n \Phi(-1/\tau )}$
using $T(-\mi u) \sim \mi/u = 2n \mi /\tau $
then (\ref{gAS}) (the regime of \S~3.1.1 holds as $|T| \gg n$).
Add that $\Phi(-1/\tau ) \sim -\tau $ ($\tau \approx 0$), and the integrand in (\ref{OUF}) begs for
the \emph{Laplace method as in} \S~3.1.1. Applied directly on (\ref{OUF}), this gives
$\int_0^{u_0} \e ^{-\mi \Theta _n (-\mi u)} u^{-\alpha } \d u \sim 
\e ^{-\mi \Theta _n (0)} \int_0^\infty \e ^{-\Theta '_n (0) \, u} [u^{-\alpha }+O(u^{1-\alpha })] \d u
= \G (1-\alpha ) \, [n(2n+1)]^{\alpha -1} + O(n^{2(\alpha -2)})$, resulting in
\beq
{\mathcal I}_{-,n}(\alpha ) = \frac{2\pi }{\G (\alpha )} [n(2n+1)]^{\alpha -1} +  O(n^{2(\alpha -2)}) .
\eeq
Finally then, recalling (\ref{IDD}) and ${\mathcal I}_{+,n}(\alpha )=o(1)$, 
\beq
{\mathcal I}_{q,n}(\alpha ) \sim
\Bigl( -\frac{\d }{\d \alpha } \Bigr) ^q \Bigr[ \frac{1}{4 \sin \hf \pi \alpha } \, {\mathcal I}_{-,n}(\alpha ) \Bigr] _{\alpha =1} \pmod{o(1)}
\eeq
with $q=0$, resp. 1, yield 
(\ref{TH_N}), resp. (\ref{TL_N}), and consequently (\ref{LNA}). 
(All the above remainders are uniform and differentiable in $\alpha $ up to $\alpha =1$.)
\hfill $\square$

\subsection{Asymptotic or full-fledged Li's criterion for $\{\L _n \}$?}
\label{LIC}
(a heuristic parenthesis).

A full Li's crterion for the new sequence $\{\L _n \}$ would be
\[
\qquad \mbox{RH} \quad \iff \quad \L _n >0 \ \mbox{ for all }n \qquad \mbox{(unproven)};
\]
while we lack a proof encompassing \emph{all} $n$, this is not essential in our focus on RH-sensitivity:
already for $\{ \l _n \}$, only $n \gg 1$ truly counted (\S~\ref{PRH}); and likewise here, 
our criterion (\ref{LNRH})--(\ref{LRH}) entails $\L _n>0$ \emph{asymptotically} if and only if RH holds 
(fig.~9 will illustrate that on a counterexample to RH).

As for lower $n$, $\L _n>0$ will be \emph{numerically} manifest there (see \S~\ref{LND}). 

All those facts still let us \emph{conjecture} that Li's criterion fully holds 
for the sequence $\{\L _n \}$ as well (but we deem it harder to prove than for $\{\l _n \}$).

\subsection{Generalized-RH asymptotic alternative}
\label{DLF}

All previous developments carry over from $\z (x)$ to \emph{some Dirichlet L-functions}: 
those of the form
\beq
\label{LDef}
L_\chi (x) \defi \sum_{k=1}^\infty \frac{\chi (k)}{k^x} \equiv d^{-x}\sum_{k=1}^d \chi (k) \, \z (x,k/d) ,
\eeq
where $\chi $ is a special type of $d$-periodic function on ${\mathbb Z}$ ($d>1$):
a \emph{real primitive Dirichlet character} ($\chi $ then has a definite parity, even or odd) 
\cite[chaps.~5--6]{Da}; $\z (x,w)$ is the Hurwitz zeta function \cite[chap.~9 eq.~(16)]{Da}).
Such $L_\chi (x)$ have very close properties to Riemann's $\z (x):$
\smallskip

\ni \textbullet \ Functional equations: 
\cite[App. B, Thm B.7 \& eq.(46)]{A} (\emph{real} primitive~$\chi $ are \emph{Kronecker symbols},
for which the usual Gaussian-sum factor $\eps (\chi) \equiv 1$)
\beq
\label{LFE}
\xi _\chi (x) \defi (\pi /d)^{-x/2} \G \bigl( \hf (x \!+\! b) \bigr) \, L_\chi (x) \equiv \xi_\chi (1-x) , \quad 
b= \biggl\{ \matrix{ 0 \ \ (\chi \mbox{ even}) \cr 1 \ \ (\chi \mbox{ odd}) \hfill} \biggr\} .
\eeq
\smallskip

\ni \textbullet \ Explicit values at integers, in terms of Bernoulli polynomials:
by the classic identities $(\ell +1) \, \zeta (-\ell ,a) = - B_{\ell +1}(a)$ ($\ell =0,1,2,\ldots$),
(\ref{LDef}) makes $L_\chi (x)$ explicit at all $x=-\ell $; 
the functional equation (\ref{LFE}) then converts that to explicit values 
for $L_\chi $ (and $\xi _\chi $) at all \emph{even, resp. odd, positive integers} 
according to the \emph{even, resp. odd} parity of $\chi $. 
Moreover, $L_\chi (1)$ (hence $\xi _\chi (1)$) is explicitly computable also for even parity
\cite[chap.~1]{Da}\cite[eq.(10.70)]{VB}.
\smallskip

\ni \textbullet \ Generalized Riemann Hypothesis (GRH): \emph{all} zeros $\rho $ of $\xi_\chi $ have $\Re \rho = \hf$.
\smallskip

\ni \textbullet \ A Keiper--Li sequence $\{ \l _{\chi ,n} \}$ exists for $\xi _\chi $ in full similarity to 
$\{ \l _n \}$ for $\xi $ (applying \cite[Cor.~1]{BL}).
For GRH, the only change in our asymptotic alternative (\ref{ASF})--(\ref{ASC}) 
affects the constant $c$ in (\ref{ASC}) or (\ref{Kas}) due to the replacement, 
in the asymptotic counting function $\overline N(T)$ from (\ref{RvM}) used in \S~\ref{KLO}, 
of the term $\log (T/\,2\pi {\rm e})$ by $\log (Td/\,2\pi {\rm e})$ \cite[chap.~16]{Da}, which results in
\beq
\label{ASL}
\mbox{GRH true: } \quad \l ^{\rm K}_{\chi ,n} \sim \hf \log n + c_d \pmod {o(1)}, \qquad c_d = c + \hf \log d .
\eeq

Then as in \S~\ref{CNS}, we can discretize the definition of $\{ \l ^{\rm K}_{\chi ,n} \}$ 
to get an explicit sequence $\{ \L _{\chi ,n} \}$ involving
finite differences of \emph{elementary} $\log \xi _\chi $-values. 
For even $\chi $, everything stays as in \S~\ref{CNS} (where even parity is implied throughout), 
hence $\L _{\chi ,n}$ are given by (\ref{EKL}) 
with $2 \, \xi (2m)$ simply replaced by $\ttfr{1}{\xi _\chi (1)} \, \xi _\chi (2m)$.
For odd $\chi $ on the other hand, the points $x_m$ have to be relocated from $2m$ to $2m+1$, 
resulting in these main changes:
\bea
\label{G-df}
[(\ref{gdf}) \to ] \qquad g^{\rm odd} _n(x) \defi \prod_{m=1}^n \frac{x + 2m}{x-2m-1}
\si\equiv\se g^{\rm odd}(x) (-1)^n \frac{\G (\hf x + 1 + n)}{\G(\hf (3 \!-\! x) \!+\! n)}, \nonumber \\[-2pt]
g^{\rm odd} (x) \si\defi\se \frac{-\sqrt \pi \, 2^{x-1}}{\cos (\pi x/2) \, \G (x)} ,
\eea
\beq
[(\ref{Anm}) \to ] \qquad A^{\rm odd} _{nm} = 2^{-2n} (n+m,n-m,2m+1)! \, / (2m + 1) , \qquad \qquad
\eeq
\beq
\label{Fodd}
[(\ref{Fdef}) \to ] \quad F^{\rm odd} _n (x) = (-1)^n \biggl[ \frac{-1}{A^{\rm odd} _{n0}} \log x
+ \sum_{m=0}^n (-1)^m A^{\rm odd} _{nm} \log (x \!-\! 2m \!-\! 1) \biggr] ,
\eeq
\beq
\label{Lodd}
[(\ref{EKL}) \to ] \quad 
\L _{\chi ,n} \equiv (-1)^n \sum_{m=1}^n (-1)^m A^{\rm odd} _{nm} \log \biggl[ \frac{1}{\xi _\chi (1)} \xi _\chi (2m \!+\! 1) \biggr]
\ \ (\chi \mbox{ odd}) .
\eeq
All in all, the asymptotic alternative (\ref{LNRH})--(\ref{LRH}) generalizes to one depending on $\chi $, 
via its parity for GRH false, and period $d$ for GRH true:
\bea
\label{LNDH}
\bullet \mbox{ GRH false: } \L _{\chi ,n} \si\sim\se \!\! 
\sum _{\rho ' - \frac{1}{2} \, \in \, 2n \, {\mathcal D}_{R_0} } \!\! F^\# _n (\rho ') \pmod{o(R^n) \ \forall R > R_0 >1} \quad \\
\mbox{and,} \!\!\!\si\se\!\!\! \mbox{for each given GRH-violating zero } \rho ' , \nonumber\\
\label{SDH1}
F^\# _n (\rho ') \si\sim\se 
\frac{g^\# (\rho ')}{\rho '(\rho '-1)} (-1)^n \frac{n^{\rho '- \frac{1}{2}} }{\log n} \quad (n \to \infty) , \\
\mbox{with} \quad (F^\# _n , g^\# ) \si\defi\se \biggl\{ \matrix{ 
(F_n , \quad g) & \mbox{for $\chi $ even, cf. (\ref{Fdef}), (\ref{gdf})} \cr
(F ^{\rm odd} _n , g^{\rm odd} ) & \mbox{for $\chi $ odd, \ cf. (\ref{Fodd}), (\ref{G-df})} } \biggr\} ; \\[6pt]
\label{LDH}
\bullet \mbox{ GRH true: } \L _{\chi ,n} \si\sim\se \log n + C_d \pmod {o(1)}, \qquad C_d = C + \hf \log d 
\eea
($\equiv c_d + \hf \log 2$), cf. (\ref{LRH}), (\ref{ASL}).
\medskip

\ni Remark: for a \emph{Siegel zero} if any ($T=0$), (\ref{SDH1}) simply gives 
$F^\# _n (\rho ') \propto \ttfr{(-1)^n n^t}{\log n}$.

\section{Quantitative aspects}

We now discuss how the discretized Keiper sequence $\{ \L _n \}$ vs the Keiper--Li $\{ \l _n \}$ 
might serve as a \emph{practical} probe for RH in a complementary perspective to standard tests.
(To rigorously (dis)prove RH using $\{ \L _n \}$ is also a prospect in theory, but we have not explored that.)

\subsection{Numerical data in the Riemann case}
\label{LND}

\begin{figure}[h]
\vskip 4mm
\center
\includegraphics[scale=.4]{voros_fig5.eps}
\caption{\small 
The coefficients $\L _n$ computed by (\ref{EKL}) up to $n=4000$, on a logarithmic $n$-scale. 
(On all plots, the line segments in gray connect data points only to aid the eye.)
Straight line: the RH-true asymptotic form (\ref{LRH}), $(\log n +C)$.}
\end{figure}

\begin{figure}[h]
\center
\includegraphics[scale=.4]{voros_fig6.eps}
\caption{\small The remainder sequence $\delta \L _n = \L _n -(\log n+C)$, 
and a rectified form $(-1)^n \delta \L _n$ (drawn for $n>300$, shifted downwards for clarity) 
to cancel the period-2 oscillations of $\delta \L _n$ at large~$n$.}
\end{figure}

Low-$n$ calculations of $\L _n$ (fig. 5) agree very early with our predicted \emph{RH-true} behavior 
(here, (\ref{LRH})), just as they did for $\l _n$ \cite{K}\cite{M1} 
and for the same reason:
given the current value $T_0 \approx 2.4 \cdot 10^{12}$ up to which RH is verified, 
an imprint (\ref{LABS}) by an RH-violating zero $\rho '$ might be of visible size
only beyond some much higher $n$-value. 
Explicit $n$-thresholds as functions of $\rho '$ will be tackled below in \S~\ref{FRH}--\ref{UCP}
(but with no unique answer like (\ref{Lup}) for $\{ \l _n \}$).

The remainder term $\delta \L _n \defi \L _n -(\log n + C)$ looks compatible with an $o(1)$ bound (fig.~6), 
albeit less neatly than the analogous $\delta \l ^{\rm K}_n$ \cite[fig.~1]{K}\cite[fig.~6b]{M1}
(note: even Keiper \cite {K} plotted $\delta \l ^{\rm L}_n = n \ \delta \l ^{\rm K}_n $). For the record,
\bea
\label{ND}
\L _1 \si\approx\se 0.069176395771 , \ \ \L _2 \approx 0.22745427267 , \ \ \L _3 \approx 0.45671413349 ; 
\ \ (\mbox{cf. } (\ref{ExV})) \nonumber \\
n \si=\se \ \ 2000 : \quad \L _n \approx 6.815360445451163 \quad (\delta \L _n \approx -0.0017849), \nonumber\\
n \si=\se \ 10000 : \quad \L _n \approx 8.428662659671506 \quad (\delta \L _n \approx +0.0020794), \nonumber\\
n \si=\se \ 20000 : \quad \L _n \approx 9.119244876955247 \quad (\delta \L _n \approx -0.000485565), \\
n \si=\se 100000 : \quad \L _n \approx 10.729678153023 \quad (\delta \L _n \approx +0.0005097985), \nonumber\\
n \si=\se 200000 : \quad \L _n \approx 11.42244991847 \quad (\delta \L _n \approx +0.0001343834) \nonumber\\
n \si=\se 500000 : \quad \L _n \approx 12.33812102688 \quad (\delta \L _n \approx -0.0004852401) \nonumber
\eea
($n>20000$ samples: courtesy of G. Misguich \cite{Mi}, see \S~\ref{HI}).

The main oscillation in $(-1)^n \delta \L _n$ (fig. 6) must come from $F_n(\rho _1)$ in (\ref{LIR})
for the lowest Riemann zero $\rho _1 :$ the $(\log n)$-period agrees with 
$2 \pi /\Im \rho _1$ (per the asymptotic form (\ref{FAS})), as $2 \pi /14.1347 \approx 0.44 \,$.

\subsection{Putative imprints of zeros violating RH}
\label{FRH}

\begin{figure}
\center
\includegraphics[scale=.4]{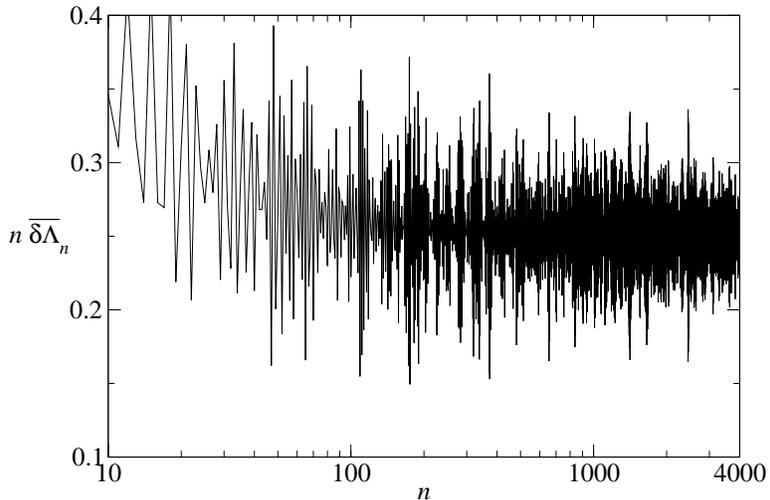}
\caption{\small The averaged remainder sequence (\ref{LAV}) rescaled by $n$, 
namely: $n \, \overline{\delta \L }_n \,$.
(Some further values: $0.27027$ for $n=10000,\ 0.23970$ for $n=20000,\ 0.2559$ for $n=100000,
\ 0.2683$ for $n=200000,\ 0.27957$ for $n=500000$ \protect\cite{Mi}.)
}
\end{figure}

RH-violating zeros $\rho '$ (if any) seem to enter the picture just as for $\{ \l _n \}$:
their contributions (\ref{ST1}) will asymptotically dominate $\log n$ from (\ref{LRH}),
but numerically they will emerge and take over extremely late. 
Indeed for such a zero $\rho ' = {\hf +t+ \mi T}$, with $0<t<\hf$ and $T \gtrsim 2.4 \cdot 10^{12}$ \cite{G},
its contribution (\ref{LABS}) scales as $T^{-2}(2n/T)^t /\log n$ in modulus. 
We then get a manifest or ``naive" crossover threshold 
(in order of magnitude, neglecting logarithms and constants against powers) by solving 
\bea
\label{TN}
T^{-2}(n/T)^t \si\approx\se 1 \\
\label{TN0}
\Longrightarrow \qquad n \si\gtrsim\se T ^{1+2/t} \qquad 
(\mbox{best case: }O(T^{5+0 }) \mbox{ for } t=\hf-0) .
\eea
This is worse than (\ref{LIt}) for $\l _n$, all the more if $\L _n <0$ were sought 
(the right-hand side of (\ref{TN}) should then be $\log ^2 n$). 
However, (\ref{TN}) is a \emph{sufficient but not necessary} condition, 
which leaves room for possible improvement.
The core problem is to filter out a weak $\rho '$-signal from the given background (\ref{LRH}),
therefore any predictable structure in the latter is liable to boost the gain.
For instance, the hyperfine structure of $\delta \L _n$ is oscillatory of period 2 (fig.~6);
this suggests to average over that period, which \emph{empirically}
discloses a rather neat $(1/n)$-decay trend (fig.~7):
\beq
\label{LAV}
\overline{\delta \L }_n \defi \hf (\delta \L _n + \delta \L _{n-1}) \approx 0.25 /n .
\eeq
The same averaging on a $\rho '$-signal $F_n (\rho ')$ in (\ref{ST1}) roughly applies $\hf (\d /\d n)$
to the factor $n^T$ therein (again neglecting $t \ll T$ and $\log n$), i.e., multiplies it by $\hf (T/n)$.
Thus heuristically, i.e., conjecturing the truth of (\ref{LAV}) for $n \to \infty$ under RH, 
the crossover condition improves from (\ref{TN}) to
\bea
\label{TNI}
(T/n) \, T^{-2}(n/T)^t \si\approx\se \overline{\delta \L }_n \approx 1 /n \nonumber \\
\Longrightarrow \quad n \si \gtrsim\se T^{1+1/t} 
\ (\mbox{best case: }O(T^{3+0 }) \mbox{ for } t=\hf-0) . \quad
\eea
We can hope that more selective signal-extraction schemes may still lower this detection threshold.
Even an empirical mindset may be acceptable there: the $\{ \L _n \}$ (like $\{ \l _n \}$) 
should anyway work best for \emph{global coarse detection} of a possible RH-violation
(as discussed in more detail in the next \S).

\subsection{The uncertainty principle for $\{ \L _n \}$}
\label{UCP}

An absolute detection limit is however provided by the \emph{uncertainty principle}, 
and this has to be quantified here. 
Near the critical line $\{ \Re X=0 \}$ parametrized by $U \defi \Im X$, and for large~$n$, 
$\L _n$ is an integral of the zeros' distribution against essentially 
$\e ^{\mi n \Psi (U)}$ (by (\ref{Lfed}), (\ref{FAN}), (\ref{PsD})): 
i.e., a high-frequency \emph{distorted plane wave}.
Then, to resolve $t>0$ for a zero $\rho ' = {\hf + t + \mi T} \equiv \hf +2nX$, 
the uncertainty principle asks its scale ($\Re X = t/(2n)$ in the $X$-plane) 
to exceed the inverse of the local wavenumber for $\e ^{\mi n \Psi (U)}$ which is 
$n \, \d \Psi /\d U$, or ${n \log (1 +1/U^2)}$ using (\ref{MOM}), with $ U= \Im X = T/(2n) :$
that means ${\hf t \, \log (1+4n^2/T^2)} \gtrsim 1$ or, to a fair approximation when $t<\hf $,
\beq
\label{UPL}
n \gtrsim \hf \, |T| \e^{1/t} .
\eeq
Against the corresponding bound (\ref{Lup}) for $\{ \l _n \}$ ($n \gtrsim T^2/t$),
(\ref{UPL}) will at any given $t>0$ favor $\{ \L _n \}$ once $|T| \gtrsim \hf t \e^{1/t}$.
E.g., at the current floor height~(\ref{TV}), the best possible $n$-threshold (i.e., for $t=\hf -0$) 
gets improved to $10^{13}$ (\emph{possibly reachable?}), 
from $10^{25}$ for $\{ \l _n \}$ in (\ref{LIt}).
The bound (\ref{UPL}) also allows $n$ to go well below
the ``naive" detection threshold (\ref{TN0}) and even (\ref{TNI}) - \emph{in principle}:
(\ref{UPL}) is a \emph{necessary but not sufficient} condition,
and there still remains to actually see the $\rho '$-signal, very weak at such decreased $n$, 
amidst the ``noise" $\delta \L_n $ brought into (\ref{LIR}) by the many more nearby zeros 
that are \emph{on} the critical line.

That $\rho '$-signal, precisely (\ref{LABS}) plus its complex-conjugate (from ${\rho '}^\ast $), is $(-1)^n \times $
a \emph{growing} oscillation of frequency $T/\, 2\pi $ in the variable ($\log n$)
(such as plotted on fig.~9, a later counterexample to $\Re \rho \equiv \hf$).
From $n$ high enough as in (\ref{TN0}) (e.g., $\gtrsim 100$ in fig.~9) this signal will overwhelmingly grow, 
and a single $\L _n$ will (most likely) betray it as evidence of an RH-violating zero \emph{somewhere},
with only (\ref{TN0}) (read backwards) to narrow its localization. 
Detection \emph{by shape-matching} with (\ref{LABS}) should be more precise, 
thus permit lower $n$-values (where the $\rho '$-signal is weaker) at the cost of using a \emph{range} of~$n$.
Better knowledge of the structure and true size of the $o(1)$ ``noise" in (\ref{LRH}) 
would be a key factor to neater results.
Barring that, we only risk highly intuitive and tentative guesswork for a \emph{buried} $\rho '$-signal. 
The $n$-variable is Fourier-conjugate to the $\Psi$-coordinate on the critical line $\{ \Re X=0 \}$; 
then the Heisenberg inequality (again) roughly bounds the localization (uncertainty) 
on $\Psi $ as $\Delta \Psi \gtrsim 1/\Delta n$ 
where $\Delta n$ is the $n$-interval on which the signal ($\L _n$) is probed;
computing as above reduces that to $\Delta U \gtrsim \log (1 +1/U^2) / \Delta n$
or, to a fair approximation under~(\ref{UPL}),
\beq
\Delta T \gtrsim 4 \log \frac{2n}{T} \, \frac{n}{\Delta n} \quad 
\Bigl( \gtrsim \frac{2T}{t} \, \frac{\e^{1/t}}{\Delta n} \mbox{ by (\ref{UPL})} \Bigr) .
\eeq
The best localization of $T$ thus improves as the width of the $n$-region grows.

All that is clearly \emph{complementary} to the direct (Riemann--Siegel based) verification methods,
which are \emph{local} in $T$ \cite{G}. These are fully rigorous, computationally more efficient, 
but strictly confined to the $T$-region on which they are applied
(thus their actual cost of use may still skyrocket if the region to be searched is unbounded).
Like $\{ \l _n \}$ but more easily, $\{ \L _n \}$ offers a ``dual" global viewpoint:
it could lead \emph{in fewer steps} to a suspicion of an RH-violating zero (if any) but neither rigorously 
nor very sharply, nevertheless enough to hint at a smaller region, in which a direct algorithm \cite{G} 
could then faster and \emph{in full rigor} confirm (or disprove) that violation of RH.

\subsection{The Davenport--Heilbronn counterexamples}
\label{DHF}

More tests of interest require the generalization of \S~\ref{DLF} to \emph{odd parity}. 
Among Dirichlet L-functions we only tested the $\beta $-function (of period~4)
and saw no difference of behavior patterns in $\{ \L _{\beta ,n} \}$ vs $\{ \L _n \}$.
If we then go beyond, there exist special Dirichlet series that are \emph{not} Dirichlet L-functions,
but obey similar functional equations \emph{and} have many zeros \emph{off} the critical line~$L$ 
\cite{DaH}\cite[\S~10.25]{Ti}\cite{S}\cite{BS}\cite[\S~5]{BG}.
In those we may seek a testing ground for the \emph{RH-false branch} of our asymptotic alternative
(generalized, as~(\ref{LNDH})--(\ref{SDH1})).

Specifically, for $\phi = \hf (1+\sqrt 5)$ (the golden ratio), let
\beq
\tau_\pm \defi -\phi \pm \sqrt{1+\phi ^2} \quad
(\tau _+ \approx 0.2840790438, \ \ \tau _- = -1/\tau _+ \approx -3.5201470213) ;
\eeq
\[
\nu _\pm (k) \defi \{ 1, \tau _\pm, -\tau _\pm, -1,0, \ldots \}_{k=1,2,\ldots} 
\quad \mbox{periodically continued}
\] 
(an \emph{odd} function on the integers mod 5); and, similarly to (\ref{LDef}),
\beq
f_\pm (x) \defi \sum_{k=1}^\infty \frac{\nu _\pm (k)}{k^x} \equiv \textstyle 
5^{-x} \{ \zeta(x,\of) + \tau_\pm \, [ \zeta(x,\tf)- \zeta(x,\frac{3}{5}) ] - \zeta(x,\frac{4}{5}) \} .
\eeq

\begin{figure}
\center
\includegraphics[scale=.4]{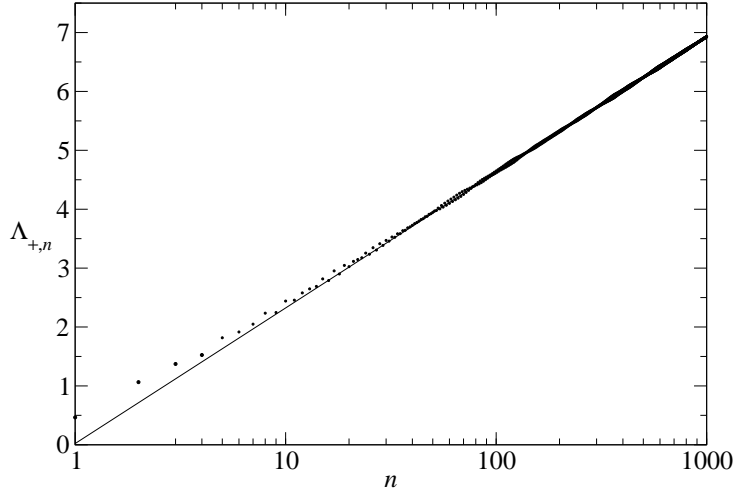}
\caption{\small As fig.~5, but for the $\L _{+,n}$ of the Davenport--Heilbronn (DH) function $f_+$
as given by (\ref{DHL}), up to $n=1000$; e.g., $\L _{+,1} \approx 0.4653858106$, $\L _{+,2} \approx 1.063986745$.
Straight line: the GRH-true form $(\log n +C_5)$ of (\ref{DHC}).}
\end{figure}

\begin{figure}[h]
\center
\includegraphics[scale=.4]{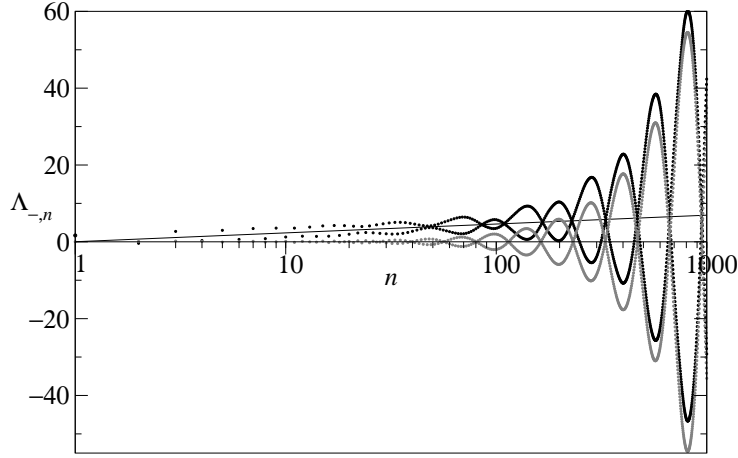}
\caption{\small As fig.~8, but for the $\L _{-,n}$ of the DH function $f_-$ (vertical scale shrunk); 
e.g., $\L _{-,1} \approx 1.661697636$, $\L _{-,2} \approx -0.3913729841$ (already $<0$). 
Straight line: the GRH-true form $(\log n +C_5)$ as in fig.~8;
gray dots: the leading GRH-false form $2 \Re F^{\rm odd}_n (\rho '_-)$ 
from (\ref{Fodd}), (\ref{LZ-}), with $\rho '_- \approx 2.30862 + 8.91836 \, \mi$.}
\end{figure}

These \emph{Davenport--Heilbronn} (DH) functions $f_\pm $ 
(denoted $f_{1 \atop 2}$ in \cite{BS}, $f(\cdot \, ; \tau_\pm)$ in \cite{BG}) 
obey \emph{the functional equation of odd Dirichlet L-functions}, 
namely (\ref{LFE}) with the parity $b=1$ and the period $d=5$, up to a ($\pm$) sign: \cite{BS}\cite[\S~5]{BG} 
\beq
\label{DHFE}
\xi _\pm (x) \defi (\pi /5)^{-x/2} \, \G [\hf(1+x)] \, f_\pm(x) \equiv \pm \, \xi_\pm (1-x) .
\eeq
As in \S~\ref{DLF}, this makes $f_\pm, \ \xi _\pm $ \emph{explicit} 
at the positive \emph{odd} integers:
\bea
\xi _\pm (2m+1) \si=\se \mp \frac{2 (-1)^m}{(2m+1)!!} \, \bigl[ B_{2m+1}(\of) + \tau _\pm B_{2m+1}(\tf) \bigr] \sqrt \pi \, (10 \, \pi)^m \, ; \\
e.g., \quad \xi _\pm (1) \si=\se \pm \of (3 + \tau_\pm) \sqrt \pi \qquad 
\Bigl\{ \matrix{\xi _+(1) \approx 1.1641757096 \cr \xi _-(1) \approx 0.1843873182} \Bigr\} \nonumber
\eea
(we also used $B_{2m+1}(1-a) \equiv -B_{2m+1}(a)$). So, in line with (\ref{Lodd}) we take
\beq
\label{DHL}
\L _{\pm,n} \defi (-1)^n \sum_{m=1}^n (-1)^m A^{\rm odd} _{mn} 
\, \log \Biggl[ \frac{1}{\xi_\pm(1)} \xi_\pm (2m \!+\! 1) \Biggr]
\eeq
as explicit sequences to test our asymptotic criterion for GRH (as (\ref{LNDH})--(\ref{LDH}) 
with $\# =$ odd and $d=5$) upon the DH functions $f_\pm$ respectively 
(as instances featuring zeros off the critical line $L$).

Now the numerics break the formal $(\pm )$-symmetry to a surprising extent.

For $\xi _+$, the lowest-$T$ zero off the line $L$ is 
$\rho '_+ \approx 0.808517 + 85.699348 \, \mi$ \cite{S}.
Then, for its detection through the sequence $\{ \L _{+,n} \}$, our predicted threshold (\ref{TN}) gives
$n \approx T^{1+2/t} \approx (85.7)^{7.48} \approx 3 \cdot 10^{14} :$ 
indeed, our low-$n$ data (fig.~8) solely reflect the GRH-true pattern (\ref{LDH}) for $d=5$,
\beq
\label{DHC}
\L _{+,n} \approx \log n + C_5 , \qquad C_5 = C + \hf \log 5 \approx +0.020961845743 \, .
\eeq

\begin{figure}
\center
\includegraphics[scale=.4]{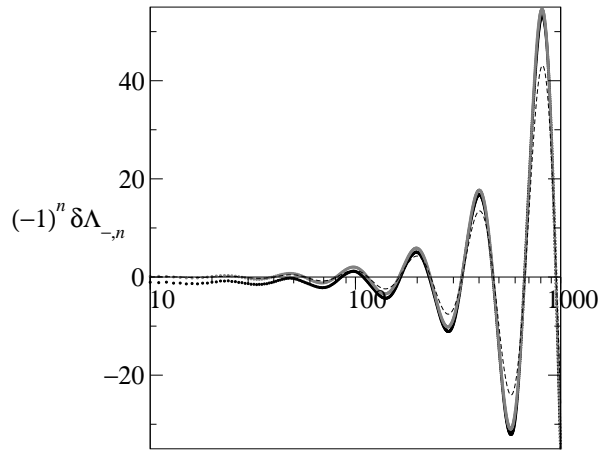}
\caption{\small Asymptotics of the deviation (\ref{DLDH}) from GRH in the case of the DH function $f_-$.
Black dots: the rectified remainder sequence $(-1)^n \delta \L _{-,n}$.
Gray dots: the rectified asymptotic form $(-1)^n \, 2 \Re F^{\rm odd}_n (\rho '_-)$ from (\ref{LZ-}). 
Dashed curve: its large-$(\log n)$ explicit form 
$2 \Re \Bigl[ \frac{g^{\rm odd} (\rho '_- )}{\rho '_- (\rho '_- -1)}n^{\rho '_- -1/2} \Bigr] (\log n)^{-1}$
from (\ref{G-df}), (\ref{SDH1}).}
\end{figure}

Whereas for $\xi _-$, the lowest-$T$ zero off $L$ is $\rho '_- \approx {2.30862 + 8.91836 \, \mi}$
\cite{BS}.
(Notations therein: $\xi \equiv \tau _+,\ f_2 \equiv f_-$.) Then, not only is this zero well detached
from the next higher one ($\approx 1.94374 + 18.8994 \,\mi$), but above all it gives a detection threshold
$n \approx T^{1+2/t} \approx (8.92)^{2.11} \approx 100$, extremely low!
Indeed, fig.~9 shows a neat crossover of $\L _{-,n}$ from the GRH-true pattern (\ref{DHC}) at low $n$,
toward the dominant GRH-false pattern (\ref{LNDH}) at higher $n$, which is 
\beq
\label{LZ-}
F^{\rm odd}_n (\rho '_-) + F^{\rm odd}_n ({\rho '}^\ast _-) \equiv 2 \Re F^{\rm odd}_n (\rho '_-) .
\eeq
Even though the form (\ref{DHC}) is asymptotically suppressed relatively to (\ref{LZ-}), 
in absolute size it remains wholly present in $\L _{-,n}$
which rather follows the straight \emph{sum} (\ref{LZ-})+(\ref{DHC}) (fig.~9).
Hence to test (\ref{LZ-}) on $\L _{-,n}$ we first have to subtract (\ref{DHC}) from the latter; 
then just like $\delta \L _n$ before (fig.~6), the remainder
\beq
\label{DLDH}
\delta \L _{-,n} \defi \L _{-,n} - (\log n + C_5)
\eeq
oscillates symmetrically about 0 with period 2, so we rather plot $(-1)^n \delta \L _{-,n}$ 
on fig.~10: we then see a very good fit by $(-1)^n \, 2 \Re F^{\rm odd}_n (\rho '_-)$ from (\ref{LZ-}).
In turn, (\ref{SDH1}) gives the latter an explicit form to first order in ($1/\log n$),
which also fits the data well enough (granted that ($1/\log n$) is not so small here).

Thus, numerical data for the sequences $\{ \L _{\pm,n} \}$ support our asymptotic alternative in full. 
We stress that fig.~9 also models how, at some much higher~$n$, $\{ \L _n \}$ itself \emph{will} blow up 
if RH (for $\z (x)$) is false.
Finally, $\{ \L _{+,n} \}$ could be a testing ground for any ideas to detect RH-violations earlier with~$\{ \L _n \}$.

\subsection{\emph{The} hitch}
\label{HI}

We now refocus on the numerical quest for $\{ \L _n \}$ (the Riemann case), to discuss a major issue.
The $(\log n)$-sized values $\L _n$ result from alternating summations like (\ref{EKL}) 
over terms growing \emph{much faster} (exponentially) with~$n$: 
very deep cancellations must then take place, which may explain the sensitivity of $\L _n$ (to RH)
but also create a computational hurdle for $n \gg 1$, namely a loss of precision growing linearly with~$n$. 
The same issue was seen on $\{ \l _n \}$ empirically \cite[fig.~6]{M2}\cite[\S~4.2]{J};
here the explicit form of $\L _n$ will give us more detailed control over that as well.

\begin{figure}[h]
\center
\vskip 5mm
\includegraphics[scale=.4]{voros_fig11.eps}
\caption{\small Base decimal precisions needed for the summands of $\L _n$ in (\ref{EKL}),
as estimated by $\log_{10} |A_{nm} \log 2\xi (2m)|$ which is plotted against $m$ in axes rescaled by $1/n$.
Dotted curve: the case $n=200$; continuous curve: the $n \to \infty$ limiting form 
$\varpi = -2 \, r \log_{10} r + (1+r) \log_{10} (1+r) -(1-r) \log_{10} (1-r) \quad (r=m/n)$.
(Cf. (\ref{PD}) and another,
seemingly unrelated, zeta-function calculation \protect\cite[\S~6.3 \& fig.~8]{BM}.)
}
\end{figure}

For sums $S=\sum s_m$ like (\ref{EKL}) that come out of order comparable to unity, 
the \emph{slightest} end accuracy requires each summand $s_m$ to be input 
with a ``base" precision $\approx \log_{10} |s_m|$ (working in decimal digits throughout); 
to which $D$ must be added (uniformly) to aim at $\approx D$ digits of accuracy for~$S$.

We can then tune the required precision in (\ref{EKL}) for each $m$-value 
at large given $n$ by using the Stirling formula,
to find that $m_\ast \approx n/\sqrt 2$ is where $|s_m|$ is largest 
and the required base precision $\log_{10} |s_m|$ culminates, reaching
$\log_{10} |A_{nm_\ast } \log 2 \xi (2m_\ast )|$ $\sim \log_{10} (3+2 \sqrt 2) \, n 
\approx 0.76555 \, n$ digits, see fig.~11 
(i.e., $2.5431 \, n$ bits, vs $n$ bits for $\l _n$ \cite[\S~4.2]{J}).
Even then, a crude feed of (\ref{EKL}) (or (\ref{Udef}), (\ref{Vdef}))
into a mainstream arbitrary-precision system (Mathematica 10 \cite{W}) 
suffices to reap the $\L _n$-values of \S~\ref{LND} effortlessly under $n \approx 20000$.
For instance (aiming at $D=13$ digits), the following returns $\L _{10000}$ 
(but would cater for $n$ \emph{of any size in principle}, by rewriting the first line only):
\smallskip

\ni \texttt{\small In[1]:= n= 10000 ; pr= 7656 + 13 ; \\
In[2]:= Timing[ \ Sum[(-1)\^{}m N[ Multinomial[n+m,n-m,2m]/(2m-1) \\
\phantom{In[2]:=} (Log[ Abs[ BernoulliB[2m]]]-Log[(2m-3)!!]), pr], \{m,n\}] /(-4)\^{}n \\
\phantom{In[2]:=} + N[ (1-(-2)\^{}n n!/(2n-1)!!)~Log[2 Pi]/2, pr] \ ] \\
Out[2]= \{172.295, 8.428662659671506\} }
\smallskip

\ni coding (\ref{Udef}), with a working precision of (7656+13) digits as just argued; 
all output digits ($16>D=13$) appear to be correct, and CPU time was 172~s.
(We used an Intel Xeon E5-2670 0 @ 2.6~GHz processor.)
Other samples of computing times we could reach (with some variance between sessions): 
ca 4~min for $\L _{10000}$, 43~min for $\L _{20000}$ using (\ref{Vdef}) 
with targeted precision $D=16$ digits and working precision dynamically adjusted according to fig.~11.

For higher $n$-values, G. Misguich kindly developed a much faster parallel code (available on request),
based on the multiple-precision GNU MPFR library \cite{Mi},
and ran it on a 20-core machine still @ 2.6 GHz with 256~Go of memory. 
He noted that the Bernoulli numbers caused the major part of the workload.
He reached CPU times $\approx $ 97~s for $n=20000$, 5.6~h for $n=100000$, 22.7 days for $n=500000$ 
(without optimizing the precision as in fig.~11; thus $n$ could still be raised,
but not so far as to justify the much longer programming and computing times then needed).

Now the true current challenge is to probe $|T|\gtrsim 2.4 \cdot 10^{12}$ by (\ref{TV}), 
hence to reach $n \gtrsim 2 \cdot 10^{36}$
(assuming the more favorable estimate (\ref{TNI}), $10^{60}$ otherwise), 
which then needs a working precision $\gtrsim 1.6 \cdot 10^{36}$ decimal places at times.
This need for a huge precision already burdened the original $\l _n$ 
but somewhat less and amidst several steeper complexities;
now for the $\L _n$, the ill-conditioning increased while other difficulties waned.
As current status, the $n$-range needed for new tests of RH stays beyond reach for the $\L _n$ too.

Still, $\{\L _n\}$ has assets to win over $\{\l _n\}$ in the long run. The $\L _n$ are fully explicit;
their evaluations are not recursive in $n$, 
thus \emph{very few} samples (at high enough~$n$, for sure)
might suffice to signal that RH is violated \emph{somewhere};
and the required working precision peaking at $\approx 0.766 \, n$ stands as the \emph{only} stumbling block,
but this is a \emph{purely logistic} barrier, 
which might be lowered if (\ref{EKL}) ever grew better conditioned variants.
Already in (\ref{Udef}), $\log 2\pi $ needs much less precision ($\propto \hf \log_{10} n$) 
as its coefficient $(2A_{n0})^{-1} \sim -\hf \sqrt{\pi n}$ grows negligibly, 
compared to the ${A_{nm} \log(|B_{2m}|/(2m-3)!! :}$ thus only the latter \emph{simpler} expressions demand 
top precision, and mainly for $m \approx n/ \sqrt 2$ (fig.~11); 
on average over $m$, that demand drops to $(2 \log_{10} 2) \, n \approx 0.602 \ n$ digits 
($\equiv 2n$ bits, vs $n$ for $\{ \l _n\}$ \cite{J}). And, to ease the computational burden of the $B_{2m} :$
a)~these may be computed just once, to evaluate many $\L _n$ in a row; 
b)~only (\ref{BE}) for $2 \xi (2m)$ uses $B_{2m}$ but \emph{numerically} speaking, for large~$m$ 
(\ref{elr})~becomes a much better alternative with its series more and more truncated as~$m$ grows. 
Further improvements of that sort should be pursued in priority to make much higher~$n$ accessible. 
The slow growth of the uncertainty-principle threshold (\ref{UPL}) 
is also an aspect favoring~$\{ \L _n \}$.

\subsection*{Concluding remark and acknowledgments}

While other sequences sensitive to RH for large $n$ are known \cite{BD}\cite{FV}, not to mention Keiper--Li again,
we are unaware of any previous case combining a fully \emph{closed form} like (\ref{EKL})
with a practical sensitivity-threshold of \emph{tempered growth} $n=O(T^\nu )$.
\medskip

We are very grateful to G. Misguich (from our Institute) who wrote and ran a special fast code 
for numerical calculations reaching $n=500000$ \cite{Mi};
and to the Referee, for many stimulating comments and suggestions which have helped us to improve this text.

\begin{appendix}
\section*{Appendix: Centered variant}

We sketch a treatment parallel to the main text for our Li-type sequences 
using the alternative basepoint $x_0=\hf$, the \emph{center} for the $\xi $-function \cite[\S 3.4]{VK}
(and focusing on the Riemann zeros' case, just for the sake of definiteness).

We recall that the Functional Equation $\xi (1-x) \equiv \xi (x)$ allows us,
in place of the mapping $z \mapsto x=(1-z)^{-1}$ within $\xi $ as in (\ref{KEd}), to use 
the double-valued map $y \mapsto x_{\tilde w} (y)= \hf \pm \sqrt{\tilde w} \, y^{1/2} / (1-y)$ 
(parametrized by $\tilde w>0$) on the unit disk.
That still maps the unit circle ${\{ |y|=1 \}}$ to the completed critical line $L \cup \{\infty\}$, 
but now minus its interval $\{ |\Im x|< \hf \sqrt{\tilde w} \}$.
As before, we ask all Riemann zeros on $L$ to pull back to ${\{ |y|=1 \}}$,
which imposes $\tilde w < \tilde w_0 \defi 4 \min_\rho |\Im \rho |^2 \approx 799.1618$. 
We thus define the parametric sequence $\{ \l ^0_n (\tilde w) \}$, for $0<\tilde w<\tilde w_0$, by
\beq
\label{Cdef}
\log 2\xi \Biggl( \hf \pm \frac{\sqrt{\tilde w} \, y^{1/2}}{1-y} \Biggr) \equiv 
\log 2\xi (\hf ) + \sum_{n=1}^\infty \frac{\l ^0_n (\tilde w)}{n} \, y^n 
\eeq
(\cite[\S 3.4]{VK} where only the case $\tilde w=1 $ was detailed, \cite{Se3}); equivalently,
\beq
\label{CRs}
\frac{\l ^0_n (\tilde w)}{n} \equiv \frac{1}{2 \pi\mi} \oint \frac{\d y}{y^{n+1}} 
\log 2\xi \bigl( x_{\tilde w} (y) \bigr) , \qquad n=1,2,\ldots 
\eeq

We now build an \emph{explicit} variant for this sequence (\ref{CRs}),
similar to $\{ \L _n \}$ for $\{ \l _n^{\rm K} \}$.
First, the deformations of (\ref{CRs}) analogous to those in \S~\ref{CNS} read as
\beq
\frac{1}{2 \pi\mi} \oint \frac{\d y}{H_{y_0}(y) \cdots H_{y_n}(y)} 
\log 2\xi (x) \qquad (\mbox{here } x \equiv x_{\tilde w} (y)) ,
\eeq
cf. (\ref{MT}), for which the simplest analytical form we found, mirroring (\ref{KLD}), is 
\beq
\label{CR}
\frac{1}{2 \pi\mi} \oint \frac{2 \, \d r}{(r+1)^2} \prod_{m=0}^n \frac{r+r_m}{r-r_m} 
\, \log 2\xi (x) , \qquad r_m \defi \frac{1+y_m}{1-y_m} ,
\eeq
where now $x \equiv x(r)$ with the new variable 
\beq
r \defi \frac{1+y}{1-y} \equiv \bigl[ 1+(2x-1)^2/{\tilde w} \bigr] ^{1/2} \qquad (\Re r>0) .
\eeq
Then with $x_m \equiv 2m$ as before (but now including $m=0$),
the integral (\ref{CR}) evaluated by the residue theorem yields the result (akin to (\ref{EKL}))
\beq
\L ^0_n (\tilde w) \defi \sum_{m=1}^n \frac{2}{(r_m \!+\! 1)^2}
\frac{\prod\limits_{k=0}^n (r_m \!+\! r_k)}{\prod\limits_{k \ne m} (r_m\!-\! r_k)} \log 2\xi (2m) , 
\quad r_m \equiv \sqrt{ 1+(4m \!-\! 1)^2/{\tilde w}} .
\eeq
These coefficients, while still explicit, are less tractable than the $\L _n$ from (\ref{EKL}), (\ref{Anm});
at the same time their design is more refined in that it captures the Functional Equation 
(through the $\l _n^0$, and unlike the $\L _n$ and $\l _n$);
nevertheless we are yet to see any \emph{concrete} benefit to using $\{ \L ^0_n \}$ over $\{ \L _n \}$.

Numerical samples (for $\tilde w = 1$, closest case to $\{ \L _n \}$; compare with (\ref{ND})):
\[
\L ^0_1 (1) \approx 0.0881535583 , \quad \L ^0_2 (1) \approx 0.237357366 , \ \ldots, 
\ \L ^0_{2000} (1) \approx 6.815307167\, .
\]

The corresponding \emph{asymptotic alternative} for RH analogous to (\ref{LNRH})--(\ref{LRH}) 
(with $\Delta_{\rho '} \L_n^0 $ meaning the contribution to $\L_n^0 $ from the zero $\rho '$)
reads as
\bea
\label{CNRH}
\bullet \ \mbox{RH false:} \ \L ^0_n (\tilde w) \si\sim\se 
\Biggl\{ \sum_{\Re \rho '>1/2} \Delta_{\rho '} \L_n^0 (\tilde w) \Biggr\}
\pmod {o(n^\eps)\ \forall \eps >0} \\
\label{ST0}
\si\se \mbox{with } \log |\Delta_{\rho '} \L_n^0 (\tilde w)| \sim (\rho ' -1/2) \log n , \\
\label{CRH}
\bullet \ \mbox{RH true:} \ \L ^0_n (\tilde w) \si\sim\se \sqrt{\tilde w} \, (\log n \!+\! C) , \quad 
C = \hf(\g \!-\! \log \pi \!-\! 1) \mbox{ as in (\ref{LRH})} . \qquad
\eea
The latter is proved by extending Oesterl\'e's method just as with $\L _n$; whereas
the former needs large-$n$ estimations of the product in (\ref{CR}), 
but our current ones remain crude compared to the full Stirling formula available for (\ref{G_N});
that precludes us from reaching the absolute scales of the $\Delta_{\rho '} \L_n^0 (\tilde w)$
and hence the values of $n$ from which any such terms might become detectable.

Numerically though (all our tests of \S~4 admit centered versions), 
the deviations we observed from the non-centered data (main text) were all slight, 
especially so for $n \gg 1$, aside from the overall factor $\sqrt{\tilde w}$ in (\ref{CRH}).
\end{appendix}

\end{document}